\documentclass[a4paper,10pt]{article}
\usepackage{stmaryrd}
\usepackage{amsfonts}
\usepackage{bbm}
\usepackage{mathrsfs}
\usepackage{latexsym,amssymb,amsmath,amscd,amsthm,amsxtra}
\usepackage[dvips]{graphicx}
\usepackage[utf8]{inputenc}
\usepackage[T1]{fontenc}
\usepackage{lmodern}
\usepackage{xcolor}
\usepackage{color}
\definecolor{red}{rgb}{1,0,0}
\definecolor{green}{rgb}{0,1,0}
\definecolor{blue}{rgb}{0,0,1}
\definecolor{refkey}{gray}{.625}
\definecolor{labelkey}{gray}{.625}
\usepackage[all]{xy}
\usepackage{nicefrac,mathtools,enumitem}
\usepackage{microtype}
\newtheorem{thm}{Theorem}[section]
\newtheorem{lem}[thm]{Lemma}

\newtheorem{pro}[thm]{Proposition}
\newtheorem{ex}[thm]{Example}
\newtheorem{rmk}[thm]{Remark}
\newtheorem{defi}[thm]{Definition}

\newcommand {\emptycomment}[1]{} 

\newcommand{\be }{\begin{equation}}
\newcommand{\ee }{\end{equation}}

\newcommand{\pf}{\noindent{\bf Proof.}\ }

\newcommand{\bl}{\bar{l} }

\newcommand{\huaA}{\mathcal{A}}
\newcommand{\huaL}{\mathcal{L}}

\newcommand{\huaE}{\mathcal{E}}

\newcommand{\huaM}{\mathcal{M}}

\newcommand{\huaV}{\mathcal{V}}

\newcommand{\huaD}{\mathcal{D}}

\newcommand{\CWM}{C^{\infty}(M)}

\newcommand{\CWhuaM}{C^{\infty}(\huaM)}

\newcommand{\frka}{\mathfrak a}

\newcommand{\frkd}{\mathfrak d}
\newcommand{\frke}{\mathfrak e}

\newcommand{\frkl}{\mathfrak l}

\newcommand{\frkX}{\mathfrak X}

\def\qed{\hfill ~\vrule height6pt width6pt depth0pt}

\newcommand{\half}{\frac{1}{2}}

\newcommand{\Courant}[1]{\left\llbracket  #1\right\rrbracket }
\newcommand{\Poisson}[1]{\{ #1\}}

\newcommand{\br}[1]{   [ \cdot,    \cdot  ]   }
\newcommand{\degree}[1]{   \mid    #1  \mid  }

\newcommand{\id}{\rm{id}}

\newcommand{\g}{\mathfrak g}

\newcommand{\dM}{\mathrm{d}}

\newcommand{\LWX}{\mathsf{CLWX}}

\newcommand{\ad}{\mathrm{ad}}

\newcommand{\sgn}{\mathrm{sgn}}
\newcommand{\Ksgn}{\mathrm{Ksgn}}

\begin{document}
\title{  QP-structures of degree 3 and   $\LWX$ $2$-algebroids\thanks
 {
Research supported by NSFC (11471139), NSF of Jilin Province (20170101050JC), Nanhu Scholars Program for Young Scholars and Nanhu Scholar Development Program of XYNU.
 }
 }
\author{Jiefeng Liu and Yunhe Sheng\\\vspace{2mm}
Department of Mathematics, Jilin University,
Changchun 130012, Jilin, China
\\ Email:~liujf12@126.com; shengyh@jlu.edu.cn
}

\date{}
\footnotetext{{\it{Keyword}: Courant algebroid, QP-structure, $\LWX$ $2$-algebroid, Leibniz $2$-algebra, Lie $3$-algebra, split Lie $2$-algebroid, split Lie $2$-bialgebroid}}
\maketitle
\begin{abstract}
 In this paper, we give the notion of a $\LWX$ 2-algebroid and show that a  QP-structure (symplectic NQ structure) of degree 3  gives rise to a $\LWX$ 2-algebroid. This is the higher analogue of the result that a  QP-structure of degree 2 gives rise to a Courant  algebroid. A $\LWX$ 2-algebroid can also be viewed as a categorified Courant algebroid. We show that one can obtain a Lie 3-algebra from a $\LWX$ 2-algebroid. Furthermore, $\LWX$ 2-algebroids are constructed from split Lie 2-algebroids and split Lie 2-bialgebroids.
\end{abstract}
\tableofcontents
\section{Introduction}

This paper is motivated by the following questions:
\begin{itemize}
  \item A QP-structure of degree $2$ gives rise to a Courant  algebroid. What is the geometric structure underlying a  QP-structure of degree 3?

  \item What is a categorified Courant algebroid? Or, equivalently, what is the $L_\infty$-analogue of a Courant algebroid?

  \item Split Lie $2$-algebroids have become a useful tool to study problems related to NQ-manifolds. What is a split Lie $2$-bialgebroid? What is the double of a split Lie $2$-bialgebroid?
\end{itemize}
The $\LWX$ 2-algebroid that we introduce in this paper provides answers of above questions.

 A QP-manifold of degree $n$ is a graded manifold equipped with a graded symplectic structure of degree $n$ and a degree $n+1$ function satisfying the master equation. A QP-manifold is also called a symplectic NQ manifold in some literature, e.g. \cite{royt}. QP-manifolds are very important in the topological field theory. Classical QP-manifolds of degree 1 are in one-to-one correspondence with Poisson manifolds. The 2-dimensional topological field theory constructed by AKSZ formulation \cite{aksz} is the Poisson sigma model. Classical QP-manifolds of degree 2 are in one-to-one correspondence with Courant algebroids \cite{royt}. Courant algebroids can be used as target spaces for a general class of 3-dimensional topological field theory \cite{RoyCF}. The notion of a Courant algebroid was introduced by Liu, Weinstein and Xu in \cite{lwx} in the study of the double of a Lie bialgebroid \cite{MackenzieX:1994}. An alternative definition was given in \cite{Roytenbergphdthesis}. See the review article \cite{Schwarzbach4} for more information. Roughly speaking, a Courant algebroid is a vector bundle, whose section space is a Leibniz algebra, together with an anchor map and a nondegenerate symmetric bilinear form, such that some compatibility conditions are satisfied. If a skew-symmetric bracket is used, in \cite{rw}, the authors showed that the underlying algebraic structure of a Courant algebroid is a Lie 2-algebra, which is the categorification of a Lie algebra \cite{baez:2algebras,Roytenberg5}. 

In \cite{Ikeda}, the authors studied QP-manifolds of degree 3 and derived a new 4-dimensional topological field theory by the AKSZ construction. The authors showed that a  QP-manifold of degree 3 gives rise to a Lie algebroid up to homotopy (Ikeda-Uchino algebroid), and analyzed its algebraic and geometric structures.

In this paper, we restudy QP-manifolds of degree 3 and find that a QP-manifold of degree 3 can give rise to a more fruitful geometric structure, which we call a $\LWX$ 2-algebroid. Roughly speaking, a $\LWX$ 2-algebroid is a graded vector bundle $\huaE=E_0\oplus E_{-1}$ over $M$, whose section space is a Leibniz 2-algebra, together with an anchor map $\rho:E_0\longrightarrow TM$ and a nondegenerate graded symmetric bilinear form of degree $1$, such that some compatibility conditions are satisfied. See Definition \ref{defi:Courant-2 algebroid} for details. Since Leibniz 2-algebras are the categorification of Leibniz algebras, $\LWX$ 2-algebroids can be viewed as the categorification of Courant algebroids. This viewpoint can also be justified by another fact: a Courant algebroid over a point is a quadratic Lie algebra while a $\LWX$ 2-algebroid over a point is a quadratic Lie 2-algebra. Generalizing Li-Bland and Meinrenken's construction of a Courant algebroid from a coisotropic action of a quadratic Lie algebra on a manifold  \cite{DLB},  we construct a $\LWX$ 2-algebroid, called the transformation $\LWX$ 2-algebroid, using an action of a quadratic Lie 2-algebra on a manifold. We show that we can obtain a Lie 3-algebra (3-term $L_\infty$-algebras) from a $\LWX$ 2-algebroid if we use the skew-symmetric bracket. This is a higher analogue of Roytenberg and Weinstein's result given in \cite{rw}.

Usually an NQ-manifold of degree $n$ is considered as a Lie $n$-algebroid \cite{Voronov:2010halgd}. In \cite{sz}, the authors defined split Lie $n$-algebroids using graded vector bundles. The equivalence between  the category of  split Lie $n$-algebroids and the category of NQ-manifolds of degree $n$ is given in \cite{BP}. The language of split Lie $n$-algebroids has slowly become a useful tool to study problems related to NQ-manifolds \cite{Jot,Jotz18}. There is a Courant algebroid structure on $A\oplus A^*$ associated to any Lie algebroid $A$. Similarly, we construct a $\LWX$ 2-algebroid structure on $\huaA\oplus \huaA^*[1]$ associated to any split Lie 2-algebroid $(\huaA=A_0 \oplus A_{-1},l_1,l_2,l_3,a)$. The notion of a Lie bialgebroid was introduced in \cite{MackenzieX:1994} as the infinitesimal object of a Poisson groupoid. Using the graded Poisson bracket on $T^*[3] E[1]$, where $E=A_0\oplus A^*_{-1}$, we introduce the notion of a split Lie 2-bialgebroid. Furthermore, we show that there is a $\LWX$ 2-algebroid structure on the double $\huaA\oplus \huaA^*[1]$ of a split Lie 2-bialgebroid $(\huaA,\huaA^*[1])$, which is a higher analogue of the fact that there is a Courant algebroid structure on the double $A\oplus A^*$ of a Lie bialgebroid $(A,A^*)$. Recently, the notion of an $L_\infty$-bialgebroid is introduced in \cite{BV}, which is a natural generalization of the Kravchenko’s notion of an $L_\infty$-bialgebra \cite{olga}. Even though the 2-term truncation of an $L_\infty$-algebroid is a split Lie 2-algebroid, the 2-term truncation of an $L_\infty$-bialgebroid is not a split Lie 2-bialgebroid. 

The theory of Courant algebroids is very rich, such as Dirac structure, Manin triples, reduction, and etc. We can go on to study these contents analogously. However, it is not the main purpose of this paper. We postpone these study in the future. In \cite{sheng}, transitive $\LWX$ 2-algebroids are studied in detail, and it is shown that a quadratic Lie 2-algebroid admits a $\LWX$-extension if and only if its first Pontryagin class, which is represented by a closed 5-form, is trivial. 

The paper is organized as follows. In Section 2, we recall QP-manifolds, Courant algebroids, Lie $n$-algebras, Leibniz 2-algebras and Lie 2-algebroids. In Section 3, we give the definition of a $\LWX$ 2-algebroid and analyze its properties. We construct ``transformation $\LWX$ 2-algebroid'' from a quadratic Lie 2-algebra action on a manifold. We show that a $\LWX$ 2-algebroid gives rise to a Lie 3-algebra (Theorem \ref{thm:Lie3}).
In Section 4, we construct a $\LWX$ 2-algebroid from a split Lie 2-algebroid directly (Theorem \ref{thm:stL2A}). In Section 5, we show that the degree 3 QP-manifold $T^*[3]A[1]$ gives rise to a $\LWX$ 2-algebroid through the derived bracket (Theorem \ref{thm:QPC2A}).
In Section 6, we give the definition of a split Lie 2-bialgebroid using the canonical graded Poisson bracket on $T^*[3]\huaA[1] $, where $\huaA=A_0\oplus A_{-1}$ is a graded vector bundle. Then we show that   the double $\huaA\oplus \huaA^*[1]$ of a split Lie 2-bialgebroid $(\huaA,\huaA^*[1])$ is a $\LWX$ 2-algebroid (Theorem \ref{thm:Lie2biC2}). 

\vspace{2mm}
 \noindent {\bf Acknowledgement:} We give our warmest thanks to Zhangju Liu, Alan Weinstein, Xiaomeng Xu and Chenchang Zhu for very useful comments and discussions. We also give our special thanks to the referee for very helpful suggestions that improve the paper.
\section{Preliminaries}

\subsection{QP-manifolds and Courant algebroids}

Recall that a graded manifold $\huaM$ is a sheaf of a graded commutative algebra over an ordinary smooth manifold $M$. The structure sheaf of $\huaM$ is locally isomorphic to a graded commutative algebra $C^\infty(U)\otimes S(V)$, where $U$ is an ordinary local chart of $M$, $S(V)$ is the polynomial algebra over $V$ and where $V:=\sum_{i\geq1}V_i$ is a graded vector space such that the dimension of $V_i$ is finite for each $i$.
\begin{defi}
A graded manifold $\huaM$ equipped with a graded symplectic structure $\omega$ of degree $n$ is called a {\bf $P$-manifold} of degree $n$.
\end{defi}
The structure sheaf $\CWhuaM$ of a $P$-manifold becomes a graded Poisson algebra.  The graded Poisson bracket is defined by
\begin{equation}
\{f,g\}=-\iota_{X_f}\iota_{X_g}\omega,
\end{equation}
where $f,g\in\CWhuaM$ and $X_f$ is the Hamiltonian vector field of $f$, i.e. $\iota_{X_f}\omega=-df$.
 We recall the basic properties of the graded Poisson bracket,
 \begin{eqnarray}
\Poisson{f,g}&=&-(-1)^{(\degree{f}-n)(\degree{g}-n)}\Poisson{g,f},\\
\Poisson{f,gh}&=&\Poisson{f,g}h+(-1)^{(\degree{f}-n)\degree{g}}g\Poisson{f,h},\\
\Poisson{f,\Poisson{g,h}}&=&\Poisson{\Poisson{f,g},h}+(-1)^{(\degree{f}-n)(\degree{g}-n)}\Poisson{g,\Poisson{f,h}},
\end{eqnarray}
where $\degree{f}$ is the degree of $f$ and $n$ is the degree of the symplectic structure. The degree of the Poisson bracket is $-n$.
\begin{defi}
Let $(\huaM,\omega)$ be a $P$-manifold of degree $n$. A function $\Theta\in\CWhuaM$ of degree $n+1$ is called a {\bf $Q$-structure}, if it is a solution of the classical master equation
\begin{equation}
\Poisson{\Theta,\Theta}=0.
\end{equation}
The triple $(\huaM,\omega,\Theta)$ is called a {\bf $QP$-manifold}.
\end{defi}

It is well-known that QP-manifolds of degree 2 are in one-to-one correspondence with Courant algebroids \cite[Theorem 4.5]{royt}.
\begin{defi}{\rm \cite{lwx}}
A {\bf Courant algebroid} is a vector  bundle $E$ together with a bundle map $\rho:E\longrightarrow TM$, a nondegenerate symmetric  bilinear form  $S$,
and an operation $\diamond:\Gamma(E)\times \Gamma(E)\longrightarrow \Gamma(E)$ such that for all $e_1,e_2,e_3\in{\Gamma(E)}$, the following axioms hold:
\begin{itemize}
\item[\rm(i)] $(\Gamma(E),\diamond)$ is a Leibniz algebra;
\item[\rm(ii)] $S(e_{1}\diamond{e_{1}}, e_2)   = \half\rho(e_2)S(e_{1},e_{1})$;
\item[\rm(iii)]
$\rho(e_{1})S(e_{2},e_{3})=S(e_{1}\diamond{e_{2}},e_{3})+S(e_{2},e_{1}\diamond{e_{3}})$.
  \end{itemize}
  \end{defi}
Given a  QP-manifold of degree $2$, the Courant algebroid structure is obtained by the derived bracket using the $Q$-structure $\Theta$ \cite{ royt}. See \cite{getzler:higher-derived,   Voronov1} for more information about higher derived brackets.

For a vector bundle $A$,   the graded manifold $T^*[2]A[1]$ is a P-manifold of degree 2.  Let $(x^i,
\xi^a)$ be local coordinates on $A[1]$, we denote by $(x^i,
\xi^a, \theta_a,p_i)$ the local coordinates on  $T^*[2]A[1]$. About their degrees, we have
$$
\mbox{degree}(x^i,
\xi^a, \theta_a,p_i)=(0,1,1,2).
$$
The graded Poisson bracket satisfies
$$
\{x^i,p_j\}=\delta^i_j=-\{p_j,x^i\},\quad \{\xi^a,\theta_b\}=\delta^a_b=\{\theta_b,\xi^a\}.
$$
A Lie algebroid structure on $A$ is equivalent to a degree $3$ function $\mu=\rho^i_bp_i\xi^b+\half\mu^a_{bc}\xi^b\xi^c\theta_a$ such that $\{\mu,\mu\}=0.$
A {\bf Lie  bialgebroid} structure on $A$ is given by a degree $3$ function $\mu+\gamma $, which can be locally written as
$$
\mu=\rho^i_bp_i\xi^b+\half\mu^a_{bc}\xi^b\xi^c\theta_a,\quad \gamma=\varrho^{ib}p_i\theta_b+\half\gamma_a^{bc}\xi^a\theta_b\theta_c,
$$
and they satisfy
$$
\{\mu+\gamma,\mu+\gamma\}=0.
$$
 On $A\oplus A^*$, there is a natural Courant algebroid structure, in which the Q-structure $\Theta$ is exactly $\mu+\gamma$.

\subsection{Lie $n$-algebras,  Leibniz 2-algebras and Lie 2-algebroids}

A  Lie $2$-algebra is a $2$-vector space $C$ equipped with
 a skew-symmetric bilinear functor,  such that the Jacobi identity is controlled by a natural isomorphism,
which satisfies the coherence law of its own.
It is well-known that a Lie $2$-algebra is  equivalent
to a 2-term $L_\infty$-algebra \cite{baez:2algebras}. $L_\infty$-algebras, also called strongly homotopy Lie algebras, were introduced  in \cite{Stasheff1}. See  \cite{LadaMartin, stasheff:introductionSHLA} for more details.

\begin{defi}
An {\bf  $L_\infty$-algebra} is a graded vector space $\mathfrak{g}=\oplus_{i\in \mathbb Z}\mathfrak{g}_{-i}$ equipped with a system $\{l_k|~1\leq k<\infty\}$
of linear maps $l_k:\wedge^k\mathfrak{g}\longrightarrow \mathfrak{g}$ with degree
$\deg(l_k)=2-k$, where the exterior powers are interpreted in the
graded sense and the following relation with Koszul sign ``{\rm Ksgn}''  is
satisfied for all $n\geq0$:
\begin{equation}\label{eq:higher-jacobi}
\sum_{i+j=n+1}(-1)^{i(j-1)}\sum_{\sigma}\sgn(\sigma)\Ksgn(\sigma)l_j
(l_i(x_{\sigma(1)},\cdots,x_{\sigma(i)}),x_{\sigma(i+1)},\cdots,x_{\sigma(n)})=0.
\end{equation}
Here the summation is taken over all $(i,n-i)$-unshuffles with
$i\geq1$.
\end{defi}

People usually refer to an $L_\infty$-algebra with $\g_{-i}=0$ for all $i\geq n$ and $i<0$ as an $n$-term $L_\infty$-algebra and we will call an $n$-term $L_\infty$-algebra a Lie $n$-algebra.

As a  model for ``Leibniz
algebras that satisfy Jacobi identity up to all higher homotopies'',
  the notion of a strongly homotopy
Leibniz algebra, or a $Lod_\infty$-algebra was given in \cite{livernet} by Livernet,
 which was further studied by Ammar, Poncin and Uchino
in \cite{ammardefiLeibnizalgebra,UchinoshL}.
  In \cite{Leibniz2al}, the authors introduced the notion of a Leibniz 2-algebra, which is the categorification of a Leibniz algebra, and prove that the category of Leibniz 2-algebras and the category of 2-term $Lod_\infty$-algebras are equivalent.

\begin{defi}\label{defi:2leibniz}
  A   {\bf Leibniz $2$-algebra} $\huaV$ consists of the following data:
\begin{itemize}
\item[$\bullet$] a complex of vector spaces $\huaV:V_{-1}\stackrel{\dM}{\longrightarrow}V_0,$

\item[$\bullet$] bilinear maps $l_2:V_{-i}\times V_{-j}\longrightarrow
V_{-i-j}$, where  $0\leq i+j\leq1$,

\item[$\bullet$] a  trilinear map $l_3:V_0\times V_0\times V_0\longrightarrow
V_{-1}$,
   \end{itemize}
   such that for all $w,x,y,z\in V_0$ and $m,n\in V_{-1}$, the following equalities are satisfied:
\begin{itemize}
\item[$\rm(a)$] $\dM l_2(x,m)=l_2(x,\dM m),$
\item[$\rm(b)$]$\dM l_2(m,x)=l_2(\dM m,x),$
\item[$\rm(c)$]$l_2(\dM m,n)=l_2(m,\dM n),$
\item[$\rm(d)$]$\dM l_3(x,y,z)=l_2(x,l_2(y,z))-l_2(l_2(x,y),z)-l_2(y,l_2(x,z)),$
\item[$\rm(e_1)$]$ l_3(x,y,\dM m)=l_2(x,l_2(y,m))-l_2(l_2(x,y),m)-l_2(y,l_2(x,m)),$
\item[$\rm(e_2)$]$ l_3(x,\dM m,y)=l_2(x,l_2(m,y))-l_2(l_2(x,m),y)-l_2(m,l_2(x,y)),$
\item[$\rm(e_3)$]$ l_3(\dM m,x,y)=l_2(m,l_2(x,y))-l_2(l_2(m,x),y)-l_2(x,l_2(m,y)),$
\item[$\rm(f)$] the Jacobiator identity:\begin{eqnarray*}
&&l_2(w,l_3(x,y,z))-l_2(x,l_3(w,y,z))+l_2(y,l_3(w,x,z))+l_2(l_3(w,x,y),z)\\
&&-l_3(l_2(w,x),y,z)-l_3(x,l_2(w,y),z)-l_3(x,y,l_2(w,z))\\
&&+l_3(w,l_2(x,y),z)+l_3(w,y,l_2(x,z))-l_3(w,x,l_2(y,z))=0.\end{eqnarray*}
   \end{itemize}
\end{defi}
We usually denote a   Leibniz 2-algebra by
$(V_{-1},V_0,\dM,l_2,l_3)$, or simply by
$\huaV$.

\begin{defi}\label{defi:Lie2algebroid}
A {\bf split Lie $2$-algebroid} is a graded vector bundle $\huaA=A_0\oplus A_{-1}$ over a manifold $M$ equipped with a bundle map (called the anchor) $a:A_0\longrightarrow TM$, and  brackets $l_i:\Gamma(\wedge^i\huaA)\longrightarrow \Gamma(\huaA)$ with degree $2-i$ for $i=1,2,3$, such that
\begin{itemize}
\item[$\rm(i)$]$(\Gamma(\huaA),l_1,l_2,l_3)$ is a Lie $2$-algebra,
\item[$\rm(ii)$]$l_2$ satisfies the Leibniz rule with respect to the anchor $a$:
$$l_2(X^0,fY)=fl_2(X^0,Y)+a(X^0)(f)Y,\quad \forall X^0\in\Gamma(A_0),f\in\CWM,Y\in\Gamma(\huaA),$$
\item[$\rm(iii)$] $l_1$ and $l_3$ are $\CWM$-linear.
 \end{itemize}
\end{defi}
 Denote a Lie 2-algebroid  by $(\huaA,l_1,l_2,l_3,a).$

\begin{rmk}
  In our definition of a Lie $n$-algebroid, the section space is an $L_\infty$-algebra. In \cite{Bruce}, the author introduced a notion of an $L_\infty$-algebroid, where the section space is a superized ($\mathbb Z_2$-graded) $L_\infty$-algebra.
\end{rmk}
\begin{lem}
Let $(\huaA,l_1,l_2,l_3,a)$ be a  Lie $2$-algebroid. Then we have
\begin{eqnarray}
\label{eq:al1}a\circ l_1&=&0,\\
\label{eq:amorphism}a(l_2(X^0,Y^0))&=&[a(X^0),a(Y^0)],\quad\forall X^0,Y^0\in\Gamma(A_0).
\end{eqnarray}
\end{lem}
\pf On one hand, for all $X^0\in\Gamma(A_0)$, $X^1\in\Gamma(A_{-1})$ and $f\in C^\infty(M)$, we have
$$
l_2(l_1(X^1),fX_0)=fl_2(l_1(X^1),X_0)+a(l_1(X^1))(f)X^0.
$$
On the other hand, since $(\Gamma(\huaA),l_1,l_2,l_3)$ is a Lie $2$-algebra, we have
$$
l_2(l_1(X^1),fX_0)=l_1(l_2(X^1,fX_0))=l_1(fl_2(X^1,X_0))=fl_1(l_2(X^1,X_0)).
$$
Therefore,  we have $a(l_1(X^1))(f)X^0=0$, which implies that \eqref{eq:al1} holds.

For all $X^0, Y^0, Z^0\in\Gamma(A_0)$ and $f\in C^\infty(M)$, by
\begin{eqnarray*}
l_2(l_2(X^0,Y^0),fZ^0)+l_2(l_2(Y^0,fZ^0), X^0)+l_2(l_2(fZ^0, X^0),Y^0)&=&-l_3(X^0,Y^0,fZ^0)\\
&=&-fl_3(X^0,Y^0,Z^0),
\end{eqnarray*}
we can deduce that \eqref{eq:amorphism} holds.
\qed

\section{$\LWX$ $2$-algebroids and Lie 3-algebras}

\subsection{$\LWX$ $2$-algebroids}
In this subsection, we   introduce the notion of a $\LWX$ 2-algebroid (named
after Courant-Liu-Weinstein-Xu) and analyze its properties.

\begin{defi}\label{defi:Courant-2 algebroid}
A {\bf $\LWX$ $2$-algebroid} is a graded vector bundle $\huaE=E_{-1}\oplus E_0$ over $M$ equipped with a non-degenerate graded symmetric bilinear form\footnote{Here graded symmetry means $S(e^i,h^j)=(-1)^{ij}S(h^j,e^i)$ for all $e^i\in\Gamma(E_{-i}),h^j\in\Gamma(E_{-j})$.} $S$ on $\huaE$, a bilinear operation $\diamond:\Gamma(E_{-i})\times \Gamma(E_{-j})\longrightarrow \Gamma(E_{-(i+j)})$, $0\leq i+j\leq 1$, which is skewsymmetric  on $\Gamma(E_0)\times \Gamma(E_0)$, an $E_{-1}$-valued $3$-form $\Omega$ on $E_0$, two bundle maps $\partial:E_{-1}\longrightarrow E_0$ and $\rho:E_0\longrightarrow TM$, such that $E_{-1}$ and $E_0$ are isotropic and the following axioms are satisfied:
\begin{itemize}
\item[$\rm(i)$]$(\Gamma(E_{-1}),\Gamma(E_0),\partial,\diamond,\Omega)$ is a Leibniz $2$-algebra;
\item[$\rm(ii)$]for all $e\in\Gamma(\huaE)$, $e\diamond e=\frac{1}{2}\huaD S(e,e)$, where $\huaD:\CWM\longrightarrow \Gamma(E_{-1})$ is defined by
 \begin{equation}
 S(\huaD f,e^0)=\rho(e^0)(f),\quad \forall e^0\in\Gamma(E_0);
 \end{equation}
\item[$\rm(iii)$]for all $e^1_1,e^1_2\in\Gamma(E_{-1})$, $S( \partial(e^1_1),e^1_2)=S(e^1_1,\partial(e^1_2))$;
\item[$\rm(iv)$]for all $e_1,e_2,e_3\in\Gamma(\huaE)$, $\rho(e_1)S( e_2,e_3)=S( e_1\diamond e_2,e_3)+S(e_2,e_1\diamond e_3)$;
\item[$\rm(v)$]for all $e^0_1,e^0_2,e^0_3,e^0_4\in\Gamma(E_0)$, $S(\Omega(e^0_1,e^0_2,e^0_3),e^0_4)=-S(e^0_3,\Omega(e^0_1,e^0_2,e^0_4))$.
 \end{itemize}
\end{defi}
 Denote a $\LWX$ 2-algebroid  by $(E_{-1},E_0,\partial,\rho,S,\diamond,\Omega)$, or simply by $\huaE$. Since the section space of a $\LWX$ 2-algebroid is a Leibniz 2-algebra, the section space of a Courant algebroid is a Leibniz algebra and  Leibniz 2-algebras are the categorification of Leibniz algebras, we can view $\LWX$ 2-algebroids   as the categorification of Courant algebroids.

 \begin{rmk}\label{rmk:a point}
   When $M$ is a point, both $E_0$ and $E_{-1}$ are vector spaces and the operators $\huaD$ and $\rho$ vanish. In this case, the operation $\diamond$ is skew-symmetric. It follows that $(E_{-1},E_0,\partial,\diamond,\Omega)$ is a Lie $2$-algebra. Furthermore, $S$ is a degree $1$ pairing. Axioms {\rm(iii)-(iv)} imply that $S$ is invariant. Thus, what we obtain is a  {\bf metric (quadratic) Lie 2-algebra}.  This is a higher analogue of the fact that a Courant algebroid over a point is a metric (quadratic) Lie algebra. See  \cite{BSZ} and \cite{olga} for more information about general notion of an $L_\infty$-algebra with a degree $k$ nondegenerate graded symmetric invariant bilinear form.
 \end{rmk}

 \begin{rmk}
   Note that via the nondegenerate bilinear form $S$, we obtain that $E_{-1}\cong E_0^*.$ Comparing to the Lie algebroid up to homotopy introduced in \cite{Ikeda},  the main difference is that our bilinear operation $\diamond$ is defined from $\Gamma(E_{-i})\times \Gamma(E_{-j})$ to $ \Gamma(E_{-(i+j)})$, $0\leq i+j\leq 1$, while their bilinear operation $[\cdot,\cdot]$ is only defined from  $\Gamma(E_0)\wedge \Gamma(E_0)$ to $\Gamma(E_0)$. Consequently, we have a Leibniz $2$-algebra underlying a $\LWX$ $2$-algebroid, which is the higher analogue of the fact that there is a Leibniz algebra underlying a Courant algebroid. It turns out that the operation $\diamond:\Gamma(E_{-i})\times \Gamma(E_{-j})\longrightarrow \Gamma(E_{-(i+j)})$, $i+j=1$, behaves more like the Courant-Dorfman bracket in a Courant algebroid. Thus, $\LWX$ $2$-algebroids are more fruitful structures than Lie algebroids up to homotopy.
 \end{rmk}

 \begin{rmk}\label{rmk:CC2}
 The   standard Courant algebroid $TM\oplus T^*M$ can be viewed as a $\LWX$ $2$-algebroid $(T^*[1]M,TM,\partial=0,\rho={\id}, S,\diamond,\Omega=0)$, where $S$ is the natural symmetric pairing between $TM$ and $T^*M$, and $\diamond$ is the standard Dorfman bracket given by
  \begin{equation}
  (X+\alpha)\diamond(Y+\beta)=[X,Y]+L_X\beta-\iota_Y d\alpha,\quad\forall~X,Y\in\frkX(M),~\alpha,\beta\in\Omega^1(M).
  \end{equation}
  Similarly, a Courant algebroid $A\oplus A^*$, in which $A$ is a Lie algebroid and $A^*$ is abelian, can also be viewed as a $\LWX$ $2$-algebroid. However, there is not a canonical way to obtain a $\LWX$ $2$-algebroid from an arbitrary Courant algebroid. See Remark \ref{rmk:QP} for an interpretation from the viewpoint of QP-manifolds.
 \end{rmk}

 \begin{ex}\label{ex:4-form}{\rm
   Let $H\in\Omega^4(M)$ be a closed $4$-form, which can be viewed as a bundle map from $\wedge^3TM\longrightarrow T^*M$. Then $(T^*[1]M,TM,\partial=0,\rho={\id}, S,\diamond,\Omega=H)$ is a $\LWX$ $2$-algebroid, where $S$ and $\diamond$ are the same as the ones given in the above remark.}
 \end{ex}

 \begin{lem}
Let $(E_{-1},E_0,\partial,\rho,S,\diamond,\Omega)$ be a $\LWX$ $2$-algebroid. For all $e_1,e_2\in\Gamma(\huaE)$, $e^0_1,e^0_2 \in\Gamma(E_0)$ and $ f\in C^\infty(M)$, we have
\begin{eqnarray}
e_1\diamond fe_2&=&f(e_1\diamond e_2)+\rho(e_1)(f)e_2,\label{eq:anchor1}\\
(fe_1)\diamond e_2&=&f(e_1\diamond e_2)-\rho(e_2)(f)e_1+ S(e_1,e_2)\huaD f,\label{eq:anchor2}\\
\rho(e^0_1\diamond e^0_2)&=&[\rho(e^0_1),\rho(e^0_2)].\label{eq:relation2}
\end{eqnarray}
\end{lem}
 \pf By axiom (iv) in   Definition \ref{defi:Courant-2 algebroid} and the nondegeneracy of $S$,  we have
 \begin{eqnarray*}
S( e_1\diamond fe_2,e_3)&=&\rho(e_1)S( fe_2,e_3)-S(fe_2,e_1\diamond e_3)\\
&=&f\rho(e_1)S(e_2,e_3)+S(e_2,e_3)\rho(e_1)(f) -fS(e_2,e_1\diamond e_3)\\
&=& S(f(e_1\diamond e_2),e_3)+S(\rho(e_1)(f)e_2,e_3),
 \end{eqnarray*}
 which implies that \eqref{eq:anchor1} holds.

By axiom (ii) in Definition \ref{defi:Courant-2 algebroid}, \eqref{eq:anchor2} follows immediately.

By (d) in Definition \ref{defi:2leibniz}, for $f\in\CWM$, we have
 \begin{eqnarray*}
f\partial\Omega(e_1^0,e_2^0,e_3^0)&=&e_1^0\diamond(e_2^0\diamond fe_3^0)-(e_1^0\diamond e_2^0)\diamond fe_3^0-e_2^0\diamond(e_1^0\diamond fe_3^0)\\
&=&f\big(e_1^0\diamond(e_2^0\diamond e_3^0)-(e_1^0\diamond e_2^0)\diamond e_3^0-e_2^0\diamond(e_1^0\diamond e_3^0)\big)\\
&&+\Big(\rho(e_1^0)\rho(e_2^0)(f)-\rho(e_2^0)\rho(e_1^0)(f)-\rho(e_1^0\diamond e_2^0)(f)\Big)e^0_3\\
&=&f\partial\Omega(e_1^0,e_2^0,e_3^0)+\Big([\rho(e^0_1),\rho(e^0_2)](f)-\rho(e_1^0\diamond e_2^0)(f)\Big)e^0_3,
 \end{eqnarray*}
which implies that \eqref{eq:relation2} holds.\qed

\begin{lem}
Let $(E_{-1},E_0,\partial,\rho,S,\diamond,\Omega)$ be a $\LWX$ $2$-algebroid. For all $e^0\in\Gamma(E_0)$ and $f\in C^\infty(M)$,  we have
\begin{eqnarray}
\rho\circ \partial&=&0,\label{eq:relation1}\\
\partial\circ\huaD&=&0,\label{eq:relation3}\\
\label{eq:DFbracketright}e^0\diamond \huaD f&=&\huaD S(e^0,\huaD f),\\
\label{eq:DFbracketleft}\huaD f\diamond e^0&=&0.
\end{eqnarray}
\end{lem}
\pf By (c) in Definition \ref{defi:2leibniz} and \eqref{eq:anchor1}, for all  $e^1_1,e^1_2\in\Gamma(E_{-1})$, we have
$$
\rho(\partial(e^1_1))(f)e^1_2=(\partial(e^1_1))\diamond (fe^1_2)-f\partial(e^1_1)\diamond e^1_2=e^1_1\diamond \partial(fe^1_2)-f\partial(e^1_1)\diamond e^1_2=0,
$$
which imply that \eqref{eq:relation1} holds.

By axiom (iii) in   Definition \ref{defi:Courant-2 algebroid} and \eqref{eq:relation1}, \eqref{eq:relation3} follows immediately.

Finally, for all $h^0\in\Gamma(E_0)$, by axiom (iv) in  Definition \ref{defi:Courant-2 algebroid} and \eqref{eq:relation2}, we have
 \begin{eqnarray*}
\rho(e^0)\rho(h^0)(f)&=&\rho(e^0)S(\huaD f,h^0)=S(e^0\diamond \huaD f,h^0)+S(\huaD f,e^0\diamond h^0)\\
&=&S(e^0\diamond \huaD f,h^0)+\rho(e^0\diamond h^0)(f)\\
&=&S(e^0\diamond \huaD f,h^0)+\rho(e^0)\rho(h^0)(f)-\rho(h^0)\rho(e^0)(f).
\end{eqnarray*}
Hence,
$$S(e^0\diamond \huaD f,h^0)=\rho(h^0)\rho(e^0)(f)=S(h^0,\huaD S(e^0,\huaD f)).$$
Since $S$ is nondegenerate, we deduce that  \eqref{eq:DFbracketright} holds.

By axiom (ii) in   Definition \ref{defi:Courant-2 algebroid}, \eqref{eq:DFbracketleft} follows immediately.\qed

\subsection{Transformation $\LWX$ 2-algebroids}

One can obtain a transformation Courant algebroid from a coisotropic action of a quadratic Lie algebra on a manifold, see \cite{DLB} for more details.
The notion of an $L_\infty$-algebra action on a graded manifold was given by Mehta and Zambon in \cite{MZ}. One can obtain a transformation $L_\infty$-algebroid from an $L_\infty$-algebra action.  Here we give explicit formulas of a Lie 2-algebra action on a usual manifold and the corresponding transformation Lie 2-algebroid, by which we construct a $\LWX$ 2-algebroid, called the transformation $\LWX$ 2-algebroid.
\begin{defi}
  An action of a Lie $2$-algebra $\g=(\g_{-1},\g_0,l_1,l_2,l_3)$ on a manifold $M$ is a linear map $\rho:\g_0\longrightarrow \frkX(M)$ such that
  \begin{eqnarray}
  \rho (l_2(x^0,y^0))&=&[\rho(x^0),\rho(y^0)],\quad \forall x^0,y^0\in\g_0,\\
  \rho\circ l_1&=&0.
  \end{eqnarray}
\end{defi}

Let $\rho:\g\longrightarrow \frkX(M)$ be an action of a Lie $2$-algebra $\g$ on a manifold $M$. Then $\rho$ induces a bundle map from $M\times\g_0$ to $TM$, which we use the same notation $\rho$. On the graded bundle $(M\times\g_{-1})\oplus (M\times \g_0)$, define
$
\bl_1:M\times\g_{-1}\longrightarrow M\times \g_0,~ \bl_2:\Gamma(M\times\g_{-i})\times \Gamma(M\times \g_{-j})\longrightarrow \Gamma(M\times\g_{-i-j}), 0\leq i+j\leq 1,$ and $ \bl_3:\wedge^3(M\times \g_0)\longrightarrow M\times \g_{-1}$ by
\begin{equation}\left\{\begin{array}{rcl}
  \bl_1(X^1)&=&l_1(X^1),\\
  \bl_2(X^0,Y^0)&=&l_2(X^0,Y^0)+L_{\rho(X^0)}Y^0-L_{\rho(Y^0)}X^0,\\
  \bl_2(X^0,Y^1)&=&-\bl_2(Y^1,X^0)=l_2(X^0,Y^1)+L_{\rho(X^0)}Y^1,\\
  \bl_3(X^0,Y^0,Z^0)&=&l_3(X^0,Y^0,Z^0).
  \end{array}\right.
\end{equation}
Then $(M\times\g_{-1},M\times\g_0,\rho,\bl_1,\bl_2,\bl_3)$ is a Lie 2-algebroid, called the transformation Lie 2-algebroid. See \cite{MZ} for the general case of transformation $L_\infty$-algebroids.

Now let $\g=(\g_{-1},\g_0,l_1,l_2,l_3)$ be a quadratic Lie $2$-algebra, i.e. there is a degree 1 nondegenerate graded symmetric invariant bilinear form $S$ on $\g$. In this case $\g_{-1}$ is isomorphic to $\g_0^*$. More precisely, the invariant condition reads
\begin{eqnarray*}
  S(l_1(x^1),y^1)&=&S(x^1,l_1(y^1)),\\
  S(l_2(x^0,y^0),z^1)&=&-S(y^0,l_2(x^0,z^1)),\\
  S(l_3(x^0,y^0,z^0),w^0)&=&-S(z^0,l_3(x^0,y^0,w^0)),
\end{eqnarray*}
for all $x^0,y^0,z^0,w^0\in\g_0$ and $x^1,y^1,z^1\in\g_{-1}$.
Let $\rho:\g\longrightarrow \frkX(M)$ be an action of   $\g$ on   $M$. With the same notations as above, on the graded bundle $(M\times\g_{-1})\oplus (M\times \g_0)$, we define the operation $\diamond:\Gamma(M\times\g_{-i})\times \Gamma(M\times \g_{-j})\longrightarrow \Gamma(M\times\g_{-i-j}), 0\leq i+j\leq 1,$ by
\begin{equation}\label{eq:fomulatran}\left\{\begin{array}{rcl}
    X^0\diamond Y^0&=&\bl_2(X^0,Y^0),\\
  X^0\diamond Y^1&=&\bl_2(X^0,Y^1)+\rho^*S(dX^0,Y^1),\\
    Y^1\diamond X^0&=&\bl_2(Y^1,X^0)+\rho^*S(dY^1,X^0), \end{array}\right.
\end{equation}
for all $X^0, Y^0\in \Gamma(M\times \g_0)$ and $Y^1\in\Gamma(M\times \g_{-1})$.

\begin{thm}
 Let $\g=(\g_{-1},\g_0,l_1,l_2,l_3)$ be a quadratic Lie $2$-algebra with a degree $1$ nondegenerate graded symmetric invariant  bilinear form $S$ on $\g$ and $\rho:\g_0\longrightarrow TM$ an action of $\g$ on $M$ such that
 \begin{equation}\label{eq:coisoaction}
   l_1\circ \rho^*=0,
 \end{equation}
 where $\rho^*:T^*M\longrightarrow M\times \g_{-1}$ is defined by $$S(\rho^*(\alpha),X^0)=\langle\alpha,\rho(X^0)\rangle,\quad\forall X^0\in\Gamma(M\times \g_0), \alpha\in\Omega^1(M).$$
 Then $(M\times\g_{-1},M\times \g_0,\partial=\bl_1,\rho,S,\diamond,\Omega=\bl_3)$ is a $\LWX$ $2$-algebroid, where $\diamond$ is given by \eqref{eq:fomulatran}.
\end{thm}

We call this $\LWX$ $2$-algebroid the {\bf transformation $\LWX$ $2$-algebroid.}

\pf Obviously, for all $X^0\in\Gamma(M\times \g_0)$ and $Y^1\in\Gamma(M\times \g_1)$, we have
$$
X^0\diamond Y^1+Y^1\diamond X^0=\rho^*(S(dX^0,Y^1)+S(X^0,dY^1))=\rho^*dS(X^0,Y^1),
$$
which implies that axiom (ii) in Definition \ref{defi:Courant-2 algebroid} holds.

For all $X^0, Y^0\in\Gamma(M\times \g_0)$ and $Z^1\in\Gamma(M\times \g_1)$, since $S$ is an invariant bilinear form on $\g$, we have
\begin{eqnarray*}
  S(X^0\diamond Y^0,Z^1)+S(Y^0,X^0\diamond Z^1)&=&S(l_2(X^0,Y^0)+L_{\rho(X^0)}Y^0-L_{\rho(Y^0)}X^0,Z^1)\\
  &&+S(Y^0,l_2(X^0,Z^1)+L_{\rho(X^0)}Z^1+\rho^*S(dX^0,Z^1))\\
  &=&S(L_{\rho(X^0)}Y^0,Z^1)+S(Y^0,L_{\rho(X^0)}Z^1)\\
  &=&\rho(X^0)S(Y^0,Z^1),
\end{eqnarray*}
which implies that axiom (iv) in Definition \ref{defi:Courant-2 algebroid} holds.

Also by the fact that $S$ is an invariant bilinear form on $\g$,  axioms (iii)  and (v) in Definition \ref{defi:Courant-2 algebroid} hold naturally.

Finally, we show that $(\Gamma(M\times\g_{-1}),\Gamma(M\times \g_0),\partial=\bl_1,\diamond,\Omega=\bl_3)$ is a Leibniz 2-algebra. By \eqref{eq:coisoaction}, we have
$$
\partial (X^0\diamond X^1)=\bl_1 (\bl_2(X^0,X^1)+\rho^*S(dX^0,X^1))=\bl_1 (\bl_2(X^0,X^1))=\bl_2(X^0,\bl_1 (X^1))=X^0\diamond \partial (X^1),
$$
which implies that Condition (a) in Definition \ref{defi:2leibniz} holds. Similarly, we can deduce that Condition (b) holds. Since $S$ is an invariant bilinear form on $\g$, we have
$$
\partial(X^1)\diamond Y^1=\bl_2(\bl_1(X^1),Y^1)+\rho^*S(d\bl_1(X^1),Y^1)=\bl_2(X^1,\bl_1(Y^1))+\rho^*S(dX^1,\bl_1(Y^1))=X^1\diamond \partial (Y^1),
$$
which implies that Condition (c) in Definition \ref{defi:2leibniz} holds.

Since for all $X^0,Y^0\in\Gamma(M\times \g_0)$, we have $X^0\diamond Y^0=\bl_2(X^0,Y^0)$. Thus, Condition (d) in Definition \ref{defi:2leibniz} holds naturally.

For all $X^0,Y^0\in\Gamma(M\times \g_0)$ and $Z^1\in\Gamma(M\times \g_{-1})$, by axiom (iv) in Definition \ref{defi:Courant-2 algebroid} that we have proved above, we have
\begin{eqnarray*}
  &&S\Big(X^0\diamond (Y^0\diamond Z^1)-(X^0\diamond Y^0)\diamond Z^1-Y^0\diamond (X^0\diamond Z^1)-\Omega(X^0,Y^0,\partial (Z^1)),Z^0\Big)\\
  &=&S\Big(X^0\diamond (\bl_2(Y^0, Z^1)+\rho^*S(dY^0,Z^1))-\bl_2(X^0,Y^0)\diamond Z^1\\
  &&-Y^0\diamond (\bl_2(X^0, Z^1)+\rho^*S(dX^0,Z^1))-\bl_3(X^0,Y^0,\bl_1 (Z^1)),Z^0\Big)\\
  &=&S\Big(\bl_2(X^0, \bl_2(Y^0, Z^1))+\rho^*S(dX^0,\bl_2(Y^0, Z^1))+X^0\diamond \rho^*S(dY^0,Z^1)\\
  &&-\bl_2(\bl_2(X^0,Y^0), Z^1)-\rho^*S(d\bl_2(X^0,Y^0), Z^1)-\bl_2(Y^0,\bl_2(X^0, Z^1))\\
  &&-\rho^*S(dY^0,\bl_2(X^0, Z^1))-Y^0\diamond \rho^*S(dX^0,Z^1)-\bl_3(X^0,Y^0,\bl_1 (Z^1)),Z^0\Big)\\
  &=&S\Big(\rho^*S(dX^0,\bl_2(Y^0, Z^1))+X^0\diamond \rho^*S(dY^0,Z^1)-\rho^*S(d\bl_2(X^0,Y^0), Z^1)\\
  &&-\rho^*S(dY^0,\bl_2(X^0, Z^1))-Y^0\diamond \rho^*S(dX^0,Z^1),Z^0\Big)\\
  &=&S(L_{\rho(Z^0)}X^0,Y^0\diamond Z^1-\rho^*S(dY^0,Z^1))+\rho(X^0)S(L_{\rho(Z^0)}Y^0,Z^1)-S(L_{[\rho(X^0),\rho(Z^0)]}Y^0,Z^1)\\
  &&-S(L_{\rho(Z^0)}\bl_2(X^0,Y^0), Z^1)-S(L_{\rho(Z^0)}Y^0, X^0\diamond Z^1-\rho^*S(dX^0,Z^1))-\rho(Y^0)S(L_{\rho(Z^0)}X^0,Z^1)\\
  &&+S(L_{[\rho(Y^0),\rho(Z^0)]}X^0,Z^1)\\
  &=&-S(\bl_2(Y^0,L_{\rho(Z^0)}X^0),Z^1)-S(L_{\rho(L_{\rho(Z^0)}X^0)}Y^0,Z^1)-S(L_{[\rho(X^0),\rho(Z^0)]}Y^0,Z^1)\\
  &&-S(L_{\rho(Z^0)}\bl_2(X^0,Y^0), Z^1)+S(\bl_2(X^0,L_{\rho(Z^0)}Y^0),Z^1)+S(L_{\rho(L_{\rho(Z^0)}Y^0)}X^0,Z^1)\\
  &&+S(L_{[\rho(Y^0),\rho(Z^0)]}X^0,Z^1)\\
  &=&-S\Big(L_{\rho(Z^0)}\bl_2(X^0,Y^0)-\bl_2(L_{\rho(Z^0)}X^0, Y^0)-\bl_2(X^0,L_{\rho(Z^0)}Y^0)+L_{\rho(L_{\rho(Z^0)}X^0)}Y^0\\
  &&-L_{\rho(L_{\rho(Z^0)}Y^0)}X^0+L_{[\rho(X^0),\rho(Z^0)]}Y^0-L_{[\rho(Y^0),\rho(Z^0)]}X^0,Z^1\Big)\\
  &=&0.
\end{eqnarray*}
The last equality is due to the following Lemma \ref{lem:techlem}. Thus, Condition ${\rm(e_1)}$ in Definition \ref{defi:2leibniz} holds. Similarly we can show that Conditions ${\rm(e_2)}$,  ${\rm(e_3)}$ and (f) in Definition \ref{defi:2leibniz} hold. Thus, $(\Gamma(M\times\g_{-1}),\Gamma(M\times \g_0),\partial=\bl_1,\diamond,\Omega=\bl_3)$ is a Leibniz 2-algebra. The proof is finished.
\qed

\begin{lem}\label{lem:techlem}
  For all $Z\in\frkX(M)$ and $X,Y\in\Gamma(M\times\g_0)$, we have
  \begin{equation}
    Z\bl_2(X,Y)-\bl_2(L_{Z}X, Y)-\bl_2(X,L_{Z}Y)+L_{\rho(L_{Z}X)}Y-L_{\rho(L_{Z}Y)}X+L_{[\rho(X),Z]}Y-L_{[\rho(Y),Z]}X=0.
  \end{equation}
\end{lem}
\pf If $X,Y\in\g$ are constant sections, it is obvious that the above equality holds. Generally, since $\Gamma(M\times\g_0)=\CWM\otimes \g_0$, we can assume that $X=fu, Y=gv$, where $u,v\in\g_0$ are constant sections and $f,g\in\CWM$, then it is straightforward to deduce the above equality. \qed

\subsection{Lie 3-algebras}
In this subsection we prove that we can obtain a Lie 3-algebra from a $\LWX$ 2-algebroid via skewsymmetrization.

We introduce a skew-symmetric bracket on $\Gamma(\huaE)$,
\begin{equation}\label{eq:skewbracket}
\Courant{e_1,e_2}=\half(e_1\diamond e_2-e_2\diamond e_1),\quad \forall~e_1,e_2\in\Gamma(\huaE),
\end{equation}
which is the skew-symmetrization of $\diamond$. By axiom (ii) in  Definition \ref{defi:Courant-2 algebroid}, \eqref{eq:skewbracket} can be written by
\begin{equation}\label{eq:skewbracket1}
\Courant{e_1,e_2}=e_1\diamond e_2-\half\huaD S(e_1,e_2).
\end{equation}

\begin{lem}
Let $(E_{-1},E_0,\partial,\rho,S,\diamond,\Omega)$ be a $\LWX$ $2$-algebroid. For all $e^0\in\Gamma(E_0),e^1,e^1_1,e^1_2\in\Gamma(E_{-1})$ and $f\in C^\infty(M)$, we have
\begin{eqnarray}
\label{eq:C2Arelation5}\partial\Courant{e^0,e^1}&=&\Courant{e^0,\partial(e^1)},\\
\label{eq:C2Arelation1}\Courant{\partial(e^1_1),e^1_2}&=&\Courant{e^1_1,\partial(e^1_2)},\\
\label{eq:C2Arelation3}\Courant{e^0,\huaD f}&=&\half\huaD S(e^0,\huaD f).
\end{eqnarray}
\end{lem}
\pf By (a) in  Definition \ref{defi:2leibniz} and \eqref{eq:relation3}, we have
$$
\partial\Courant{e^0,e^1}=\partial(e^0\diamond e^1)-\half\partial\circ\huaD S(e_0,e_1)=e^0\diamond \partial(e^1),
$$
which implies that \eqref{eq:C2Arelation5} holds.

By (c) in  Definition \ref{defi:2leibniz} and axiom (iii) in   Definition \ref{defi:Courant-2 algebroid}, \eqref{eq:C2Arelation1} follows immediately.

By \eqref{eq:DFbracketright} and \eqref{eq:DFbracketleft}, \eqref{eq:C2Arelation3} is obvious.  \qed\vspace{3mm}

For simplicity, for all $e_i\in\Gamma(\huaE),i=1,2,3,$ we let
\begin{eqnarray}
K(e_1,e_2,e_3)&=&e_1\diamond(e_2\diamond e_3)-(e_1\diamond e_2)\diamond e_3-e_2\diamond(e_1\diamond e_3),\\
J(e_1,e_2,e_3)&=&\Courant{\Courant{e_1,e_2},e_3}+\Courant{\Courant{e_2,e_3},e_1}+\Courant{\Courant{e_3,e_1},e_2}.
\end{eqnarray}
By \eqref{eq:relation2} and \eqref{eq:DFbracketleft}, we can deduce that $K$ is totally skew-symmetric.

\begin{lem}
Let $(E_{-1},E_0,\partial,\rho,S,\diamond,\Omega)$ be a $\LWX$ $2$-algebroid. For all $e^0,e^0_1,e^0_2,e^0_3\in\Gamma(E_0),$ $e^1, e^1_1,e^1_2\in\Gamma(E_{-1})$, we have
\begin{eqnarray}
J(e^0_1,e^0_2,e^0_3)&=&-\partial\Omega(e^0_1,e^0_2,e^0_3),\label{eq:C2AJ1}\\
J(e^0_1,e^0_2,e^1)&=&\huaD T(e^0_1,e^0_2,e^1)-\Omega(e^0_1,e^0_2,\partial e^1),\label{eq:C2AJ2}\\
T(\partial e^1_1,e^0,e^1_2)&=&-T(\partial e^1_2,e^0,e^1_1),\label{eq:C2Arelation2}
\end{eqnarray}
where the totally skew-symmetric  $T:\Gamma(E_0)\times \Gamma(E_0)\times \Gamma(E_{-1})\longrightarrow C^\infty(M)$ is given by
\begin{equation}
T(e^0_1,e^0_2,e^1)=\frac{1}{6}\big(S(e^0_1,\Courant{e^0_2,e^1})+S(e^1,\Courant{e^0_1,e^0_2})+S(e^0_2,\Courant{e^1,e^0_1})\big).
\end{equation}
\end{lem}
\pf
It is obvious that  $J(e^0_1,e^0_2,e^0_3)=-K(e^0_1,e^0_2,e^0_3)$, which implies that \eqref{eq:C2AJ1} holds.

By straightforward computations, we have
\begin{eqnarray*}
K(e^0_1,e^0_2,e^1)=-J(e^0_1,e^0_2,e^1)+R(e^0_1,e^0_2,e^1),
\end{eqnarray*}
where \begin{eqnarray*}
R(e^0_1,e^0_2,e^1)&=&\half \big(\huaD S(e^0_1,\Courant{e^0_2,e^1})-\huaD S(e^1,\Courant{e^0_1,e^0_2})
-\huaD S(e^0_2,\Courant{e^0_1,e^1})\\&&+\huaD S(e^0_1,\huaD S(e^1,e^0_2))-\huaD S(e^0_2,\huaD S(e^1,e^0_1))\big).
\end{eqnarray*}
Similarly, we have
\begin{eqnarray*}
K(e^1,e^0_1,e^0_2)=-J(e^1,e^0_1,e^0_2)+R(e^1,e^0_1,e^0_2),
\end{eqnarray*}
where \begin{eqnarray*}
R(e^1,e^0_1,e^0_2)=\half \big(\huaD S(e^0_2,\Courant{e^0_1,e^1})+\huaD S(e^1,\Courant{e^0_1,e^0_2})+\huaD S(e^0_1,\Courant{e^0_2,e^1})-\huaD S(e^0_1,\huaD S(e^1,e^0_2))\big),
\end{eqnarray*}
and
\begin{eqnarray*}
K(e^0_2,e^1,e^0_1)=-J(e^0_2,e^1,e^0_1)+R(e^0_2,e^1,e^0_1),
\end{eqnarray*}
where \begin{eqnarray*}
R(e^0_2,e^1,e^0_1)=\half \big(-\huaD S(e^0_2,\Courant{e^0_1,e^1})-\huaD S(e^0_1,\Courant{e^0_2,e^1})+\huaD S(e^1,\Courant{e^0_1,e^0_2})+\huaD S(e^0_2,\huaD S(e^1,e^0_1))\big).
\end{eqnarray*}
Since both $J$ and $K$ are completely skew-symmetric,  we have
 \begin{eqnarray*}
3K(e^0_1,e^0_2,e^1)&=&-3J(e^0_1,e^0_2,e^1)+3\huaD T(e^0_1,e^0_2,e^1).
\end{eqnarray*}
Then by axiom ${\rm(e_1)}$ in Definition \ref{defi:2leibniz}, we have
$$
K(e^0_1,e^0_2,e^1)= \Omega(e^0_1,e^0_2,\partial e^1),
$$
which implies that \eqref{eq:C2AJ2} holds.

Finally, by axiom (iii) in the Definition \ref{defi:Courant-2 algebroid}, \eqref{eq:C2Arelation5} and \eqref{eq:C2Arelation1}, we have
\begin{eqnarray*}
T(\partial (e^1_1),e^0,e^1_2)&=&\frac{1}{6}\big(S(\partial (e^1_1),\Courant{e^0,e^1_2})+S(e^0,\Courant{e^1_2,\partial (e^1_1)})+S(e^1_2,\Courant{\partial (e^1_1),e^0})\big)\\
&=&\frac{1}{6}\big(S(e^1_1,\Courant{e^0,\partial (e^1_2)})+S(e^0,\Courant{\partial (e^1_2),e^1_1})+S(\partial (e^1_2),\Courant{e^1_1,e^0})\big)\\
&=&-T(\partial e^1_2,e^0,e^1_1).
\end{eqnarray*}
The proof is finished.\qed

\begin{lem}\label{lem:proof4}
For all  $e^1\in\Gamma(E_{-1})$ and $e^0_1, e^0_2,e^0_3,e^0_4\in\Gamma(E_0)$, we have
\begin{eqnarray*}
&&\Omega(\Courant{e^0_1,e^0_2},e^0_3,e^0_4)-\Omega(\Courant{e^0_1,e^0_3},e^0_2,e^0_4)+\Omega(\Courant{e^0_1,e^0_4},e^0_2,e^0_3)
+\Omega(\Courant{e^0_2,e^0_3},e^0_1,e^0_4)\\&&-\Omega(\Courant{e^0_2,e^0_4},e^0_1,e^0_3)+\Omega(\Courant{e^0_3,e^0_4},e^0_1,e^0_2)
-\Courant{\Omega(e^0_1,e^0_2,e^0_3),e^0_4}-\Courant{\Omega(e^0_1,e^0_3,e^0_4),e^0_2}\\&&+\Courant{{\Omega}(e^0_1,e^0_2,e^0_4),e^0_3}
+\Courant{\Omega(e^0_2,e^0_3,e^0_4),e^0_1}+\huaD S(\Omega(e^0_1,e^0_2,e^0_3),e^0_4)=0,
\end{eqnarray*}
and
$$2{\bf J}+{\bf K}=- S(\Omega(\partial e^1,e^0_2,e^0_3),e^0_4),$$
where
\begin{eqnarray*}
{\bf J}&=&S(J(e^1,e^0_2,e^0_3),e^0_4)-S(J(e^1,e^0_2,e^0_4),e^0_3)+S(J(e^1,e^0_3,e^0_4),e^0_2)+3 S(\Omega(\partial e^1,e^0_2,e^0_3),e^0_4),\\
{\bf K}&=&S(\Courant{e^1,e^0_2},\Courant{e^0_3,e^0_4})-S(\Courant{e^1,e^0_3},\Courant{e^0_2,e^0_4})+S(\Courant{e^1,e^0_4},\Courant{e^0_2,e^0_3}).
\end{eqnarray*}
\end{lem}
\pf By  axiom  (f) in Definition \ref{defi:2leibniz},  axiom (v) in  Definition \ref{defi:Courant-2 algebroid} and \eqref{eq:skewbracket1}, we have
\begin{eqnarray*}
&&\Omega(\Courant{e^0_1,e^0_2},e^0_3,e^0_4)-\Omega(\Courant{e^0_1,e^0_3},e^0_2,e^0_4)+\Omega(\Courant{e^0_1,e^0_4},e^0_2,e^0_3)
+\Omega(\Courant{e^0_2,e^0_3},e^0_1,e^0_4)\\&&-\Omega(\Courant{e^0_2,e^0_4},e^0_1,e^0_3)+\Omega(\Courant{e^0_3,e^0_4},e^0_1,e^0_2)
-\Courant{\Omega(e^0_1,e^0_2,e^0_3),e^0_4}-\Courant{\Omega(e^0_1,e^0_3,e^0_4),e^0_2}\\&&+\Courant{{\Omega}(e^0_1,e^0_2,e^0_4),e^0_3}
+\Courant{\Omega(e^0_2,e^0_3,e^0_4),e^0_1}+\huaD S(\Omega(e^0_1,e^0_2,e^0_3),e^0_4)\\
&=&\Omega({e^0_1\diamond e^0_2},e^0_3,e^0_4)-\Omega({e^0_1\diamond e^0_3},e^0_2,e^0_4)+\Omega({e^0_1\diamond e^0_4},e^0_2,e^0_3)+\Omega({e^0_2\diamond e^0_3},e^0_1,e^0_4)\\
&&-\Omega({e^0_2\diamond e^0_4},e^0_1,e^0_3)+\Omega({e^0_3\diamond e^0_4},e^0_1,e^0_2)-\Omega(e^0_1,e^0_2,e^0_3)\diamond e^0_4+\half\huaD S(\Omega(e^0_1,e^0_2,e^0_3),e^0_4)\\&&+e^0_2\diamond\Omega(e^0_1,e^0_3,e^0_4)-\half \huaD S(\Omega(e^0_1,e^0_3,e^0_4),e^0_2)-e^0_3\diamond\Omega(e^0_1,e^0_2,e^0_4)+\half \huaD S(\Omega(e^0_1,e^0_2,e^0_4),e^0_3)\\
&&-e^0_1\diamond\Omega(e^0_2,e^0_3,e^0_4)+\half\huaD S(\Omega(e^0_2,e^0_3,e^0_4),e^0_1)+\huaD S(\Omega(e^0_1,e^0_2,e^0_3),e^0_4)\\
&=&\half\huaD S(\Omega(e^0_1,e^0_2,e^0_3),e^0_4)-\half \huaD S(\Omega(e^0_1,e^0_3,e^0_4),e^0_2)+\half \huaD S(\Omega(e^0_1,e^0_2,e^0_4),e^0_3)\\
&&+\half\huaD S(\Omega(e^0_2,e^0_3,e^0_4),e^0_1)+\huaD S(\Omega(e^0_1,e^0_2,e^0_3),e^0_4)\\
&=&-\huaD S(\Omega(e^0_1,e^0_2,e^0_3),e^0_4)+\huaD S(\Omega(e^0_1,e^0_2,e^0_3),e^0_4)=0.
\end{eqnarray*}
The second equality can be proved by the same method in the proof of Lemma 2.5.2 in \cite{Roytenbergphdthesis}. We omit the details.\qed\vspace{3mm}
\emptycomment{
We will extend the $\LWX$ $2$-algebroid structure to an Lie-$3$ algebra on the total space of the following resolution of $H=\rm{coker}\huaD:$
\begin{equation}\label{eq:sequence}
\cdots 0\longrightarrow X_3 \stackrel{d_3}\longrightarrow X_2\stackrel{d_2}\longrightarrow X_1\stackrel{d_1}\longrightarrow X_0\longrightarrow H\longrightarrow 0,
\end{equation}
where $X_0=\Gamma(E_0),X_1=\Gamma(E_{-1}),X_2=\CWM,X_3=\rm{ker}\huaD,d_1=\partial,d_2=\huaD$ and $d_3$ is the inclusion $\imath:\rm{ker}\huaD\hookrightarrow\CWM$.
}

Let $(E_{-1},E_0,\partial,\rho,S,\diamond,\Omega)$ be a $\LWX$ 2-algebroid. Consider the graded vector space $\frke=\frke_{-2}\oplus \frke_{-1}\oplus \frke_0$, where $\frke_0=\Gamma(E_0)$, $\frke_{-1}=\Gamma(E_{-1})$ and $\frke_{-2}=C^\infty(M)$.

\begin{thm}\label{thm:Lie3}
A $\LWX$ $2$-algebroid  $(E_{-1},E_0,\partial,\rho,S,\diamond,\Omega)$ gives rise to a Lie $3$-algebra $(\frke,l_1,l_2,l_3,l_4)$, where $l_i$ are given by the following formulas:
$$\begin{array}{rcll}
l_1(f)&=&\huaD(f),& \forall~ f\in C^\infty(M),\\
l_1(e^1)&=&\partial (e^1),& \forall~ e^1\in \Gamma(E_{-1}),\\
l_2(e^0_1\wedge e^0_2)&=&\Courant{e^0_1,e^0_2},& \forall~ e^0_1, e^0_2\in\Gamma(E_0),\\
l_2(e^0\wedge e^1)&=&\Courant{e^0,e^1},& \forall~ e^0  \in\Gamma(E_0), e^1\in\Gamma(E_{-1}), \\
l_2(e^0\wedge f)&=&\half S(e^0,\huaD f),& \forall~ e^0  \in\Gamma(E_0), f\in C^\infty(M),\\
l_2(e^1_1\vee e^1_2)&=&0,& \forall~  e^1_1, e^1_2\in\Gamma(E_{-1}),\\
l_3(e^0_1\wedge e^0_2\wedge e^0_3)&=&\Omega(e^0_1, e^0_2, e^0_3),& \forall~ e^0_1, e^0_2, e^0_3\in\Gamma(E_0),\\
 l_3(e^0_1\wedge e^0_2\wedge e^1)&=&-T(e^0_1, e^0_2, e^1),& \forall~ e^0_1, e^0_2\in\Gamma(E_0), e^1\in \Gamma(E_{-1}),\\
l_4(e^0_1\wedge e^0_2\wedge e^0_3\wedge e^0_4)&=&\overline{\Omega}(e^0_1, e^0_2, e^0_3,e^0_4),&\forall~ e^0_1, e^0_2, e^0_3, e^0_4\in\Gamma(E_0),
\end{array}
$$
where $\overline{\Omega}:\wedge^4\Gamma(E_0) \longrightarrow C^\infty(M)$ is given by
$$
\overline{\Omega}(e^0_1, e^0_2, e^0_3,e^0_4)=S(\Omega(e^0_1, e^0_2, e^0_3),e^0_4).
$$
\end{thm}
\pf We need to show that \eqref{eq:higher-jacobi} holds for $n=1,2,3,4,5$. For $n=1$, we need to show that $l_1^2=0$, which follows from  $\partial\circ\huaD=0$.

For $n=2$, we need to verify the following equality:
\begin{equation}\label{eq:Jn2}
-l_2(l_1(x_1),x_2)+(-1)^{\degree{x_1}\degree{x_2}}l_2(l_1(x_2),x_1)+l_1l_2(x_1,x_2)=0,\quad \forall~x_i\in\frke.
\end{equation}
For  $x_1=e^0\in\frke_0, x_2=f\in\frke_{-2}$, by \eqref{eq:C2Arelation3}, we have
\begin{eqnarray*}
l_2(\huaD f,e^0)+\huaD l_2(e^0, f)&=&-\Courant{e^0,\huaD f}+\half\huaD S(e^0,\huaD f)\\
&=&-\half\huaD S(e^0,\huaD f)+\half\huaD S(e^0,\huaD f)=0,
\end{eqnarray*}
which implies that \eqref{eq:Jn2} holds for $x_1\in\frke_0$ and $x_2\in\frke_{-2}$. The other cases  can be proved similarly and we omit the details.

For $n=3$, we need to prove that for all $x_i\in\frke$, the following equality holds:
\begin{eqnarray}
\nonumber&&l_3(l_1(x_1),x_2,x_3)-(-1)^{\degree{x_1}\degree{x_2}}l_3(l_1(x_2),x_1,x_3)+(-1)^{\degree{x_3}(\degree{x_1}+\degree{x_2})}l_3(l_1(x_3),x_1,x_2)\\
&&+l_2(l_2(x_1,x_2),x_3)-(-1)^{\degree{x_2}\degree{x_3}}l_2(l_2(x_1,x_3),x_2)+(-1)^{\degree{x_1}(\degree{x_2}+\degree{x_3})}l_2(l_2(x_2,x_3),x_1)\nonumber\\
\label{eq:JI3}&&+l_1l_3(x_1,x_2,x_3)=0.
\end{eqnarray}
  By \eqref{eq:C2AJ1}, we can deduce that \eqref{eq:JI3} holds for $x_1,x_2,x_3\in\frke_0.$ By \eqref{eq:C2AJ2}, we can deduce that \eqref{eq:JI3} holds for two elements in $ \frke_0$ and one element in $\frke_{-1}$.
 By \eqref{eq:C2Arelation2}, we can deduce that \eqref{eq:JI3} holds for one element in $ \frke_0$ and two elements in $\frke_{-1}.$
 The other cases   can be proved similarly and  we omit the details.

For $n=4$, we need to verify the following equality:
\begin{eqnarray*}
&&-l_4(l_1(x_1),x_2,x_3,x_4)+(-1)^{\degree{x_1}\degree{x_2}}l_4(l_1(x_2),x_1,x_3,x_4)-(-1)^{\degree{x_3}(\degree{x_1}+\degree{x_2})}l_4(l_1(x_3),x_1,x_2,x_4)\\
&&+(-1)^{\degree{x_4}(\degree{x_1}+\degree{x_2}+\degree{x_3})}l_4(l_1(x_4),x_1,x_2,x_3)-(-1)^{\degree{x_2}\degree{x_3}}l_3(l_2(x_1,x_3),x_2,x_4)\\
&&+(-1)^{\degree{x_4}(\degree{x_3}+\degree{x_2})}l_3(l_2(x_1,x_4),x_2,x_3)+(-1)^{\degree{x_1}(\degree{x_3}+\degree{x_2})}l_3(l_2(x_2,x_3),x_1,x_4)\\
&&-(-1)^{\degree{x_1}(\degree{x_4}+\degree{x_2})+\degree{x_3}\degree{x_4}}l_3(l_2(x_2,x_4),x_1,x_3)+(-1)^{(\degree{x_3}+\degree{x_4})(\degree{x_1}+\degree{x_2})}l_3(l_2(x_3,x_4),x_1,x_2)\\
&&-l_2(l_3(x_1,x_2,x_3),x_4)-(-1)^{\degree{x_2}(\degree{x_3}+\degree{x_4})}l_2(l_3(x_1,x_3,x_4),x_2)+(-1)^{\degree{x_3}\degree{x_4}}l_2(l_3(x_1,x_2,x_4),x_3)\\
&&+(-1)^{\degree{x_1}(\degree{x_2}+\degree{x_3}+\degree{x_4})}l_2(l_3(x_2,x_3,x_4),x_1)+l_3(l_2(x_1,x_2),x_3,x_4)+l_1l_4(x_1,x_2,x_3,x_4)=0.
\end{eqnarray*}
For $x_1=e^0_1,x_2=e^0_2,x_3=e^0_3,x_4=e^0_4\in\frke_0$, we need to prove that
\begin{eqnarray*}
&&\Omega(\Courant{e^0_1,e^0_2},e^0_3,e^0_4)-\Omega(\Courant{e^0_1,e^0_3},e^0_2,e^0_4)+\Omega(\Courant{e^0_1,e^0_4},e^0_2,e^0_3)\\
&&+\Omega(\Courant{e^0_2,e^0_3},e^0_1,e^0_4)-\Omega(\Courant{e^0_2,e^0_4},e^0_1,e^0_3)+\Omega(\Courant{e^0_3,e^0_4},e^0_1,e^0_2)\\
&&-\Courant{\Omega(e^0_1,e^0_2,e^0_3),e^0_4}-\Courant{\Omega(e^0_1,e^0_3,e^0_4),e^0_2}+\Courant{{\Omega}(e^0_1,e^0_2,e^0_4),e^0_3}\\
&&+\Courant{\Omega(e^0_2,e^0_3,e^0_4),e^0_1}+\huaD\overline{\Omega}(e^0_1,e^0_2,e^0_3,e^0_4)=0,
\end{eqnarray*}
which holds by Lemma \ref{lem:proof4}.

For $x_1=e^1\in\frke_{-1},x_2=e^0_2,x_3=e^0_3,x_4=e^0_4\in\frke_0$, we need to prove that
\begin{eqnarray*}
&&-\overline{\Omega}(\partial e^1,e^0_2,e^0_3,e^0_4)-T(\Courant{e^1,e^0_2},e^0_3,e^0_4)+T(\Courant{e^1,e^0_3},e^0_2,e^0_4)\\
&&-T(\Courant{e^1,e^0_4},e^0_2,e^0_3)-T(\Courant{e^0_2,e^0_3},e^0_1,e^0_4)+T(\Courant{e^0_2,e^0_4},e^0_1,e^0_3)\\
&&-T(\Courant{e^0_3,e^0_4},e^1,e^0_2)+\Courant{T(e^1,e^0_2,e^0_3),e^0_4}+\Courant{T(e^1,e^0_3,e^0_4),e^0_2}\\
&&-\Courant{T(e^0_2,e^0_3,e^0_4),e^1}-\Courant{T(e^1,e^0_2,e^0_4),e^0_3}=0.
\end{eqnarray*}
On one hand, by direct calculation, we have
\begin{eqnarray*}
\Courant{T(e^1,e^0_2,e^0_3),e^0_4}+\Courant{T(e^1,e^0_3,e^0_4),e^0_2}-\Courant{T(e^0_2,e^0_3,e^0_4),e^0_1}-\Courant{T(e^1,e^0_2,e^0_4),e^0_3}=-\half {\bf J}.
\end{eqnarray*}
On the other hand, we have
\begin{eqnarray*}
&&-\overline{\Omega}(\partial e^1,e^0_2,e^0_3,e^0_4)-T(\Courant{e^1,e^0_2},e^0_3,e^0_4)+T(\Courant{e^1,e^0_3},e^0_2,e^0_4)\\
&&-T(\Courant{e^1,e^0_4},e^0_2,e^0_3)-T(\Courant{e^0_2,e^0_3},e^1,e^0_4)+T(\Courant{e^0_2,e^0_4},e^1,e^0_3)\\
&&-T(\Courant{e^0_3,e^0_4},e^1,e^0_2)=-\frac{1}{6}({\bf J}+2{\bf K})-\frac{1}{3}\overline{\Omega}(\partial e^1,e^0_2,e^0_3,e^0_4).
\end{eqnarray*}
Therefore, by Lemma \ref{lem:proof4}, we prove the equality above.

Finally, we can show that \eqref{eq:higher-jacobi} holds for $n=5$. We omit the details. The proof is finished. \qed

\begin{rmk}
  In \cite{Roytenberg5}, Roytenberg showed that one can obtain a semistrict Lie $2$-algebra from a weak Lie $2$-algebra via the skew-symmetrization. For a $\LWX$ $2$-algebroid  $(E_{-1},E_0,\partial,\rho,S,\diamond,\Omega)$, the Leibniz $2$-algebra $(\Gamma(E_{-1}),\Gamma(E_0),\partial,\diamond,\Omega)$ is not necessarily a weak Lie $2$-algebra. Thus, we obtain a Lie $3$-algebra rather than a Lie $2$-algebra via the skew-symmetrization.
\end{rmk}

\begin{rmk}
  In this remark, we give a possible way to understand Theorem \ref{thm:Lie3} conceptually. In \cite{Roytenberg5}, Roytenberg introduced the notion of a weak Lie $2$-algebra and showed that via skew-symmetrization, one can obtain a Lie $2$-algebra. Assume that this result could be generalized to the higher case: one can obtain a Lie $n$-algebra from a weak Lie $n$-algebra via skew-symmetrization. Then hopefully our Leibniz $2$-algebra in a $\LWX$ $2$-algebroid can naturally be completed to a weak Lie $3$-algebra and the Lie $3$-algebra given in Theorem  \ref{thm:Lie3} is exactly its skew-symmetrization.
  \end{rmk}

\section{The $\LWX$ $2$-algebroid associated to a split Lie 2-algebroid}
In this section, we first  describe a split Lie 2-algebroid structure on a graded vector bundle $A_{-1}\oplus A_0$ using the graded Poisson bracket on $T^*[3](A_0\oplus A^*_{-1})[1]$. Then we construct a $\LWX$ 2-algebroid $\huaA\oplus\huaA^*[1]$ from a split Lie 2-algebroid $\huaA$ with explicit formulas using the usual language of differential calculus. In Section 6, we will generalize this result to the case of split Lie 2-bialgebroids using the tool of derived brackets and graded geometry.

Let $\huaA= A_{-1}\oplus A_0$ be a graded bundle. The shifted cotangent bundle $T^*[3](A_0\oplus A^*_{-1})[1]$  is a $P$-manifold of degree $3$ over $M$.
 Denote by $(x^i,\xi^j,\theta_k,p_i,\xi_j,\theta^k)$ a canonical Darboux coordinate on $\huaM$, where $x^i$ is a coordinate on $M$, $(\xi^j,\theta_k)$ is the fiber coordinate on $A_0\oplus A^*_{-1}$, $(p_i,\xi_j,\theta^k)$ is the momentum coordinate on $\huaM$ for $(x^i,\xi^j,\theta_k)$. The degrees of variables $(x^i,\xi^j,\theta_k,p_i,\xi_j,\theta^k)$ are respectively $(0,1,1,3,2,2)$.
 The degree of the symplectic structure $\omega=dx^idp_i+d\xi^j d\xi_j+d\theta_kd\theta^k$ is $3$ and the degree  of the corresponding graded Poisson structure is $-3$.

Now we consider the following function\footnote{We thank very much the referee for pointing out that such a function is linear on $\huaA^*$.} $\mu$ of degree $4$ on $\huaM$:
\begin{equation} \label{eq:L2Amu}
\mu={\mu_1}^i_j(x)p_i\xi^j+{\mu_2}^{i}_j(x)\xi_i\theta^j+\half{\mu_3}^k_{ij}(x)\xi_k\xi^i\xi^j+ {\mu_4}_{ij}^k\theta^j\xi^i\theta_k+\frac{1}{6}{\mu_5}_{ijk}^l(x)\theta_l\xi^i\xi^j\xi^k,
\end{equation}
where ${\mu_1}^i_j,{\mu_2}^{i}_j,{\mu_3}^k_{ij},{\mu_4}_{ij}^k,{\mu_5}_{ijk}^l$ are   functions on $ M$.
The function $\mu$ can be uniquely decomposed into\footnote{It is routine to check that the decomposition does not depend on the choice of local coordinates. See also \cite{Ikeda} for more details.}
 $$\mu=\mu_2+\mu_{134}+\mu_5,$$
 where  $\mu_2,\mu_{134}$ and $\mu_5$ are given by
\begin{eqnarray*}
\mu_2&=&{\mu_2}^{i}_j(x)\xi_i\theta^j,\\
\mu_{134}&=&{\mu_1}^i_j(x)p_i\xi^j+\half{\mu_3}^k_{ij}(x)\xi_k\xi^i\xi^j+{\mu_4}_{ij}^k\xi^i\theta^j\theta_k,\\
\mu_5&=&\frac{1}{6}{\mu_5}_{ijk}^l(x)\theta_l\xi^i\xi^j\xi^k.
\end{eqnarray*}


Define a bundle map $l_1:A_{-1}\longrightarrow A_0$  by
\begin{eqnarray}
l_1(X^1)&=&\Poisson{X^1,\mu_2}.\label{eq:L2Al1}
\end{eqnarray}
Define $l_2:\Gamma(A_{-i})\times\Gamma( A_{-j})\longrightarrow \Gamma(A_{-i-j}),~0\leq i+j\leq 1$ by
\begin{equation}\label{eq:L2Al21}\left\{\begin{array}{rcl}
l_2(X^0,Y^0)&=&\Poisson{Y^0,\Poisson{X^0,\mu_{134}}},\\
l_2(X^0,Y^1)&=&\Poisson{Y^1,\Poisson{X^0,\mu_{134}}},\\
l_2(Y^1, X^0)&=&-\Poisson{X^0,\Poisson{Y^1,\mu_{134}}}.\end{array}\right.
\end{equation}
Define a bundle map  $l_3:\wedge^3A_0\longrightarrow A_{-1}$ by
\begin{eqnarray}
l_3(X^0,Y^0,Z^0)&=&\Poisson{Z^0,\Poisson{Y^0,\Poisson{X^0,\mu_{5}}}},\label{eq:L2Al3}
\end{eqnarray}
where $X^0,Y^0,Z^0\in\Gamma(A_0)$ and $X^1,Y^1\in\Gamma(A_{-1})$.

Finally, define a bundle map $a:A_0\longrightarrow TM$ by
\begin{equation}\label{eq:L2Aa}
a(X^0)(f)=\Poisson{f,\Poisson{X^0,\mu_{134}}},\quad \forall~X^0\in\Gamma(A_0),f\in\CWM.
\end{equation}

\begin{thm}\label{thm:L2A-MST}
Let  $\huaA=A_{-1} \oplus A_0$ be a graded vector bundle and $\mu$ a degree $4$ function given by \eqref{eq:L2Amu}. If $\Poisson{\mu,\mu}=0$,   $(\huaA,l_1,l_2,l_3,a)$ is a split Lie $2$-algebroid, where $l_1$, $l_2$, $l_3$ and $a$ are given by  \eqref{eq:L2Al1}-\eqref{eq:L2Aa} respectively.

Conversely, if $(\huaA,l_1,l_2,l_3,a)$ is a split Lie $2$-algebroid, we have $\Poisson{\mu,\mu}=0$, where $\mu$ is given by \eqref{eq:L2Amu}, in which ${\mu_1}^i_j,{\mu_2}^{i}_j,{\mu_3}^k_{ij},{\mu_4}_{ij}^k,{\mu_5}_{ijk}^l$ are given by:
\begin{eqnarray*}
a(\xi_j)&=&{\mu_1}^i_j\frac{\partial}{\partial x^i},\quad l_1(\theta_j)={\mu_2}^{i}_j\xi_i,\\
l_2(\xi_i,\xi_j)&=&{\mu_3}^k_{ij}\xi_k,\quad l_2(\theta_j,\xi_i)= {\mu_4}^k_{ij}\theta_k,\quad l_3(\xi_i,\xi_j,\xi_k)={\mu_5}_{ijk}^l\theta_l.
\end{eqnarray*}
\end{thm}
\pf One can easily prove that $\Poisson{\mu,\mu}=0$ is equivalent to the following three identities:
\begin{eqnarray*}
\Poisson{\mu_{134},\mu_2}&=&0,\\
\frac{1}{2}\Poisson{\mu_{134},\mu_{134}}+\Poisson{\mu_2,\mu_5}&=&0,\\
\Poisson{\mu_{134},\mu_5}&=&0.
\end{eqnarray*}
 It is straightforward to deduce that Conditions (ii) and (iii) in Definition \ref{defi:Lie2algebroid}   holds.

In the following, we prove that $(\Gamma(\huaA),l_1,l_2,l_3)$ is a Lie $2$-algebra. It is easy to see that $l_2$ and $l_3$ are totally skew-symmetric. For all $X^0\in\Gamma(A_0),X^1\in\Gamma(A_{-1})$, we have
\begin{eqnarray*}
\Poisson{X^1,\Poisson{X^0,\Poisson{\mu_2,\mu_{134}}}}=-l_2(X^0,l_1(X^1))+l_1l_2(X^0,X^1)=0,
\end{eqnarray*}
which implies that $l_1l_2(X^0,X^1)=l_2(X^0,l_1(X^1))$.

For all $X^1,Y^1\in\Gamma(A_{-1})$, we have
 $$\Poisson{Y^1,\Poisson{X^1,\Poisson{\mu_2,\mu_{134}}}}=l_2(l_1(X^1),Y^1)-l_2(X^1,l_1(Y^1))=0,$$
 which implies that $l_2(l_1(X^1),Y^1)=l_2(X^1,l_1(Y^1))$.

For all $X^0,Y^0,Z^0\in\Gamma(A_0)$, by $$\Poisson{Z^0,\Poisson{Y^0,\Poisson{X^0,\frac{1}{2}\Poisson{\mu_{134},\mu_{134}}+\Poisson{\mu_2,\mu_5}}}}=0,$$
we get
 $$l_2(X^0,l_2(Y^0,Z^0))+l_2(Z^0, l_2(X^0,Y^0))+l_2(Y^0,l_2(Z^0, X^0))=l_1l_3(X^0,Y^0,Z^0).$$

For all $X^0,Y^0\in\Gamma(A_0),Z^1\in\Gamma(A_{-1})$, by $\Poisson{Z^1,\Poisson{Y^0,\Poisson{X^0,\frac{1}{2}\Poisson{\mu_{134},\mu_{134}}+\Poisson{\mu_2,\mu_5}}}}=0$,
we get $$l_2(X^0,l_2(Y^0,Z^1))+l_2(Z^1, l_2(X^0,Y^0))+l_2(Y^0,l_2(Z^1, X^0))=l_3(X^0,Y^0,l_1(Z^1)).$$

For all $X^0,Y^0,Z^0,W^0\in\Gamma(A_{0})$, by
\begin{eqnarray*}
\Poisson{W^0,\Poisson{Z^0,\Poisson{Y^0,\Poisson{X^0,\frac{1}{2}\Poisson{\mu_{134},\mu_{134}}+\Poisson{\mu_2,\mu_5}}}}}=0,
\end{eqnarray*}
we deduce that \eqref{eq:higher-jacobi} holds for $n=4$.
Therefore, $(\Gamma(\huaA),l_1,l_2,l_3)$ is a Lie $2$-algebra.

The proof of the converse part is similar as the above deduction. We omit the details. The proof is finished.   \qed\vspace{3mm}

Let $(\huaA,l_1,l_2,l_3,a)$ be a split Lie $2$-algebroid with the structure function $\mu$. Then we have a generalized Chevalley-Eilenberg complex $(\Gamma(Sym(\huaA[1])^*),\delta)$,
\emptycomment{
 $\Gamma(\huaA)$ is given by
$C^m(\huaA)=\oplus_{k+2l=m}\Gamma(\wedge^lA^*_0\wedge\odot^kA^*_{-1}).$ The Chevalley-Eilenberg operator $\delta:C^m(\huaA)\longrightarrow C^{m+1}(\huaA)$
 }
 where $\delta$ is defined by
\begin{equation}
\delta(\cdot)=\{\mu,\cdot\}.
\end{equation}
 In particular, for all $f\in\CWM,\alpha^0\in\Gamma(A^*_0),\alpha^1\in\Gamma(A_{-1}^*)$, we have
\begin{eqnarray}\label{eq:delta}
\left\{\begin{array}{rcl}
\delta (f)(X^0)&=&a(X^0)(f),\\
\delta(\alpha^0)(X^0,Y^0)&=&a(X^0)\langle\alpha^0,Y^0\rangle-a(Y^0)\langle\alpha^0,X^0\rangle-\langle\alpha^0,l_2(X^0,Y^0)\rangle,\\
\delta(\alpha^1)(X^0,Y^1)&=&a(X^0)\langle\alpha^1,Y^1\rangle-\langle\alpha^1,l_2(X^0,Y^1)\rangle,
\end{array}\right.
\end{eqnarray}
where $X^0,Y^0\in\Gamma(A_0),Y^1\in\Gamma(A_{-1})$.

 Given a split Lie 2-algebroid   $(\huaA,l_1,l_2,l_3,a),$ define $l_1^*:A^*_0\longrightarrow A_{-1}^*$ by
\begin{equation}
\langle l_1^*(\alpha^0),X^1\rangle=\langle\alpha^0,l_1(X^1)\rangle,\quad \forall \alpha^0\in\Gamma(A^*_0),Y^1\in\Gamma(A_{-1}).
\end{equation}

For all $X^0\in\Gamma(A_0)$, define $L^0_{X^0}:\Gamma(A^*_{-i})\longrightarrow \Gamma(A^*_{-i})$, $i=0,1$, by
\begin{eqnarray*}
\langle L^0_{X^0}\alpha^0,Y^0\rangle&=&\rho(X^0)\langle Y^0,\alpha^0\rangle-\langle \alpha^0,l_2(X^0,Y^0)\rangle,\quad \forall \alpha^0\in\Gamma(A^*_0),Y^0\in\Gamma(A_0),\\
\langle L^0_{X^0}\alpha^1,Y^1\rangle&=&\rho(X^0)\langle Y^1,\alpha^1\rangle-\langle \alpha^1,l_2(X^0,Y^1)\rangle,\quad \forall \alpha^1\in\Gamma(A_{-1}^*),Y^1\in\Gamma(A_{-1}).
\end{eqnarray*}

For all $X^1\in\Gamma(A_{-1})$, define $L^1_{X^1}:\Gamma(A_{-1}^*)\longrightarrow \Gamma(A^*_0)$ by
\begin{equation}
\langle L^1_{X^1}\alpha^1,Y^0\rangle=-\langle \alpha^1,l_2(X^1,Y^0)\rangle,\quad \forall \alpha^1\in\Gamma(A_{-1}^*),Y^0\in\Gamma(A_0).
\end{equation}

For all $X^0,Y^0\in\Gamma(A_0)$, define $L^3_{X^0,Y^0}:\Gamma(A_{-1}^*)\longrightarrow \Gamma(A^*_0)$ by
\begin{equation}
\langle L^3_{X^0,Y^0}\alpha^1,Z^0\rangle=-\langle \alpha^1,l_3(X^0,Y^0,Z^0)\rangle,\quad \forall \alpha^1\in\Gamma(A_{-1}^*),Z^0\in\Gamma(A_0).
\end{equation}

The following lemmas list some properties of the above operators.
\begin{lem}
 For all $X^0\in\Gamma(A_0),X^1\in\Gamma(A_{-1}),f\in\CWM,$ $\alpha^0\in\Gamma(A^*_0),\alpha^1\in\Gamma(A_{-1}^*), $ we have
 \begin{eqnarray*}
  L^0_{X^0}f\alpha^0&=&f(L^0_{X^0}\alpha^0)+a(X^0)(f)\alpha^0,\\
 L^0_{fX^0}\alpha^0&=&f(L^0_{X^0}\alpha^0)+\langle X^0,\alpha^0\rangle\delta(f),\\
 L^0_{X^0}f\alpha^1&=&f( L^0_{X^0}\alpha^1)+a(X^0)(f)\alpha^1,\\
 L^0_{fX^0}\alpha^1&=&f( L^0_{X^0}\alpha^1),\\
  L^1_{X^1}f\alpha^1&=&f(L^1_{X^1}\alpha^1),\\
   L^1_{fX^1}\alpha^1&=&f( L^1_{X^1}\alpha^1)+\langle X^1,\alpha^1\rangle\delta(f),\\
  L^0_{X^0}\alpha^0&=&\iota_{X^0}\delta \alpha^0+\delta\iota_{X^0}\alpha^0,\\
  L^1_{X^1}\alpha^1&=&\delta\iota_{X^1}\alpha^1-\iota_{X^1}\delta \alpha^1.
 \end{eqnarray*}
 \end{lem}
 \pf It is straightforward.\qed

 \begin{lem}\label{lem:LieJacobi}
 For $X^0,Y^0\in\Gamma(A_0),X^1\in\Gamma(A_{-1}),\alpha^0\in\Gamma(A^*_0),\alpha^1\in\Gamma(A_{-1}^*)$, we have
 \begin{eqnarray}
\label{eq:jacobi-eq1}L^0_{l_2(X^0, Y^0)}\alpha^0- L^0_{X^0}L^0_{Y^0}\alpha^0+L^0_{Y^0}L^0_{X^0}\alpha^0&=&-L^3_{X^0,Y^0}l^*_1\alpha^0,\\
 \label{eq:jacobi-eq2}L^0_{l_2(X^0, Y^0)}\alpha^1- L^0_{X^0}L^0_{Y^0}\alpha^1+L^0_{Y^0}L^0_{X^0}\alpha^1&=&-l^*_1L^3_{X^0,Y^0}\alpha^1,\\
 L^1_{l_2(X^1, Y^0)}\alpha^1-L^1_{X^1}L^0_{Y^0}\alpha^1+L^0_{Y^0}L^1_{X^1}\alpha^1&=&-L^3_{l_1(X^1),Y^0}\alpha^1.
 \end{eqnarray}
 \end{lem}
 \pf For all $Z^0\in\Gamma(A_0),$ we have
 \begin{eqnarray*}
 &&\langle L^0_{l_2(X^0, Y^0)}\alpha^0- L^0_{X^0}L^0_{Y^0}\alpha^0+L^0_{Y^0}L^0_{X^0}\alpha^0,Z^0\rangle\\
 &=&(a(l_2(X^0,Y^0))-a(X^0)a(Y^0)+a(Y^0)a(X^0))\langle\alpha^0,Z^0\rangle\\
 &&+\langle\alpha^0,-l_2(l_2(X^0,Y^0),Z^0)-l_2(Y^0,l_2(X^0,Z^0))+l_2(X^0,l_2(Y^0,Z^0))\rangle\\
 &=&\langle\alpha^0,l_1l_3(X^0,Y^0,Z^0)\rangle\\
 &=&\langle-L^3_{X^0,Y^0}l^*_1\alpha^0,Z^0\rangle,
 \end{eqnarray*}
 which implies that the first equality holds.  The others can be proved similarly. \qed\vspace{3mm}

Let $(\huaA,l_1,l_2,l_3,a)$ be a split Lie $2$-algebroid. Now let $E_0=A_0\oplus A^*_{-1}$, $E_{-1}=A_{-1}\oplus A^*_{0}$ and $\huaE=E_0\oplus E_{-1}$. Let $\partial:E_{-1}\longrightarrow E_0$ and $\rho:E_0\longrightarrow TM$ be bundle maps defined by
\begin{eqnarray}
\partial(X^1+\alpha^0)&=&l_1(X^1)+l^*_1(\alpha^0),\label{eq:stpartial}\\
\rho(X^0+\alpha^1)&=&a(X^0)\label{eq:stanchor}.
\end{eqnarray}

 On $\Gamma(\huaE)$, there is a natural symmetric bilinear form $(\cdot,\cdot)_+$ given by
\begin{equation} \label{eq:naturalsymform}
(X^0+\alpha^1+X^1+\alpha^0,Y^0+\beta^1+Y^1+\beta^0)_+=\langle X^0,\beta^0 \rangle+\langle Y^0,\alpha^0 \rangle+\langle X^1,\beta^1 \rangle+\langle Y^1,\alpha^1 \rangle,
\end{equation}
where $X^0,Y^0\in\Gamma(A_0),X^1,Y^1\in\Gamma(A_{-1}),\alpha^0,\beta^0\in\Gamma(A^*_0),\alpha^1,\beta^1\in\Gamma(A_{-1}^*)$.

On $\Gamma(\huaE)$, we introduce the operation $\diamond$ by
\begin{equation}\label{eq:0-bracket}\left\{\begin{array}{rcl}
(X^0+\alpha^1)\diamond(Y^0+\beta^1)&=&l_2(X^0,Y^0)+L^0_{X^0}\beta^1-L^0_{Y^0}\alpha^1,\\
(X^0+\alpha^1)\diamond(X^1+\alpha^0)&=&l_2(X^0,X^1)+L^0_{X^0}\alpha^0+\iota_{X^1}\delta(\alpha^1), \\
(X^1+\alpha^0)\diamond(X^0+\alpha^1)&=&l_2(X^1,X^0)+L^1_{X^1}\alpha^1-\iota_{X^0}\delta(\alpha^0).
\end{array}\right.
\end{equation}

An $E_{-1}$-valued $3$-form $\Omega$ is defined by
\begin{equation}\label{eq:st3-form}
\Omega(X^0+\alpha^1,Y^0+\beta^1,Z^0+\zeta^1)=l_3(X^0,Y^0,Z^0)+L^3_{X^0,Y^0}\zeta^1+L^3_{Z^0,X^0}\beta^1+L^3_{Y^0,Z^0}\alpha^1,
\end{equation}
where $X^0,Y^0,Z^0\in\Gamma(A_0),\alpha^1,\beta^1,\zeta^1\in\Gamma(A_{-1}^*).$

It is easy to see that the operator $\huaD:\CWM\longrightarrow \Gamma(E_{-1})$ is given by
\begin{equation}
\huaD(f)=\delta(f),\quad\forall ~f\in\CWM.
\end{equation}

\begin{thm}\label{thm:stL2A}
Let $(\huaA,l_1,l_2,l_3,a)$ be a split Lie $2$-algebroid. Then $(E_{-1},E_0,\partial,\rho,S,\diamond,\Omega)$ is a $\LWX$ $2$-algebroid, where $\partial$ is given by \eqref{eq:stpartial}, $\rho$ is given by \eqref{eq:stanchor}, $S$ is given by \eqref{eq:naturalsymform}, $\diamond$ is given by \eqref{eq:0-bracket} and  $\Omega$ is given by \eqref{eq:st3-form}.
\end{thm}
\pf It is easy to verify that $e\diamond e=\half \huaD(e,e)_+$ for all $e\in\Gamma(\huaE)$.

In the following, we verify that $(\Gamma(E_{-1}),\Gamma(E_0),\partial,\diamond,\Omega)$ is a Leibniz $2$-algebra.
For all $e^0=X^0+\alpha^1\in\Gamma(E_0),e^1=X^1+\alpha^0,$
we have
\begin{eqnarray*}
\partial((X^0+\alpha^1)\diamond(X^1+\alpha^0))&=&l_1l_2(X^0,X^1)+l_1^*\big(L^0_{X^0}\alpha^0+\iota_{X^1}\delta(\alpha^1)\big),\\
(X^0+\alpha^1)\diamond\partial(X^1+\alpha^0)&=&l_2(X^0,l_1(X^1))+L^0_{X^0}l_1^*(\alpha^0)-L^0_{l_1(X^1)}\alpha^1.
\end{eqnarray*}
Since $(\Gamma(\huaA),l_1,l_2,l_3)$ is a Lie $2$-algebra, we have
$$l_1l_2(X^0,X^1)=l_2(X^0,l_1(X^1)),\quad l_2(l_1(X^1), Y^1)=l_2(X^1, l_1(Y^1)).$$
Then by the fact that $a\circ l_1=0$, we get
\begin{eqnarray*}
l_1^*\big(L^0_{X^0}\alpha^0+\iota_{X^1}\delta(\alpha^1)=L^0_{X^0}l_1^*(\alpha^0)-L^0_{l_1(X^1)}\alpha^1.
\end{eqnarray*}
 Therefore we have
 \begin{equation}
  \partial(e^0\diamond e^1)=e^0\diamond\partial(e^1),
 \end{equation}
which implies that Condition (a) in Definition \ref{defi:2leibniz} holds.

Also by the fact  $a\circ l_1=0$, we have
\begin{eqnarray*}
\partial(e^1\diamond e^0)=l_1^*(\delta(e^1,e^0)_+)-\partial(e^0\diamond e^1)=-\partial(e^0\diamond e^1)
=-e^0\diamond\partial(e^1)=\partial(e^1)\diamond e^0,
\end{eqnarray*}
which implies that Condition (b) in Definition \ref{defi:2leibniz} holds.

Similarly, for all $e^1_i\in\Gamma(E_{-1}), i=1,2$, we have
\begin{equation}
  \partial(e^1_1)\diamond e^1_2=e^1_1\diamond\partial(e^1_2),
\end{equation}
which implies that Condition (c) in Definition \ref{defi:2leibniz} holds.

For all $X_i^0\in\Gamma(A_0),i=1,2,3$, it is obvious that
$$K(X_1^0,X_2^0,X_3^0)=l_1l_3(X_1^0,X_2^0,X_3^0)=\partial\Omega(X_1^0,X_2^0,X_3^0).$$
Furthermore, for all $X_i^0\in\Gamma(A_0),i=1,2$ and $\alpha^1\in\Gamma(A_{-1}^*)$, by Lemma \ref{lem:LieJacobi}, we have
\begin{eqnarray*}
K(X_1^0,X_2^0,\alpha^1)&=&-(L^0_{l_2(X_1^0, X_2^0)}\alpha^1- L^0_{X_1^0}L^0_{X_2^0}\alpha^1+L^0_{X_2^0}L^0_{X_1^0}\alpha^1)\\
&=&l_1^*L^3_{X_1^0,X_2^0}\alpha^1=\partial\Omega(X_1^0,X_2^0,\alpha^1).
\end{eqnarray*}
Therefore, for all $e^0_i\in\Gamma(E_0),i=1,2,3$, we get
\begin{equation}
K(e^0_1,e^0_2,e^0_3)=\partial\Omega(e^0_1,e^0_2,e^0_3),
\end{equation}
which implies that Condition (d) in Definition \ref{defi:2leibniz} holds.

Similarly, for all $e^0_i \in\Gamma(E_0),i=1,2$ and $e^1 \in\Gamma(E_{-1})$, we have
\begin{eqnarray*}
K(e^0_1,e^0_2,e^1)=\Omega(e^0_1,e^0_2,\partial e^1),\quad
K(e^0_1,e^1,e^0_2)=\Omega(e^0_1,\partial e^1,e^0_2),\quad
K(e^1 ,e^0_1,e^0_2)=\Omega(\partial e^1,e^0_1,e^0_2),
\end{eqnarray*}
which implies that Conditions ($e_1$)-($e_3$)  in Definition \ref{defi:2leibniz} holds.

By the coherence law that $l_3$ satisfies in the definition of a Lie 2-algebra, we can deduce that Condition (f) in Definition \ref{defi:2leibniz} also holds. We omit the details. Thus, $(\Gamma(E_{-1}),\Gamma(E_0),\partial,\diamond,\Omega)$ is a Leibniz $2$-algebra.

Finally, for all $e^1_1, e^1_2\in\Gamma(E_{-1})$, $e_1, e_2, e_3\in\Gamma(\huaE)$ and $e^0_1,e^0_2,e^0_3,e^0_4\in\Gamma(E_0)$, it is straightforward to deduce that
\begin{eqnarray*}
(\partial (e^1_1),e^1_2)_+&=&(e^1_1,\partial (e^1_2))_+, \\
\rho(e_1)(e_2,e_3)_+&=&(e_1\diamond e_2,e_3)_++(e_2,e_1\diamond e_3)_+, \\
(\Omega(e^0_1,e^0_2,e^0_3),e^0_4)_+&=&-(e^0_3,\Omega(e^0_1,e^0_2,e^0_4))_+,
\end{eqnarray*}
which implies that axioms (iii), (iv) and (v) in Definition \ref{defi:Courant-2 algebroid} hold. The proof is finished.\qed

\emptycomment{
\eqref{eq:split2} and \eqref{eq:split5} are direct. For $e^0_1=X^0_1+\alpha^1_1,e^0_2=X^0_2+\alpha^1_2\in\Gamma(E_0),e^1_1=Y^1+\beta^0\in\Gamma(E_{-1})$, on one hand, we have
$$\rho(e^0_1)(e^0_2,e^1_1)_+=\rho(X^0_1)\langle X^0_2,\beta^0\rangle+\rho(X^0_1)\langle Y^1,\alpha^1_2\rangle.$$
On the other hand, we have
\begin{eqnarray*}
 &&(e^0_1\circ e^0_2,e^1_1)_++(e^0_2,e^0_1\circ e^1_1)_+\\
 &=&(l_2(X^0_1,X^0_2)+L^0_{X^0}\alpha^1_2-L^0_{X^0_2}\alpha^1,Y^1+\beta^0)_+\\
 &&+(X^0_2+\alpha^1_2,l_2(X^0_1,Y^1)+L^0_{X^0_1}\beta^0-\iota_{Y^1}\delta(\alpha^1))_+\\
 &=&\langle l_2(X^0_1,X^0_2),\beta^0\rangle+\rho(X^0_1)\langle \alpha^1_2,Y^1\rangle-\langle \alpha^2_1,l_2(X^0_1,Y^1)\rangle-\rho(X^0_2)\langle \alpha^1_1,Y^1\rangle+\langle \alpha^1_1,l_2(X^0_2,Y^1)\rangle\\
 &&+\langle \alpha^2_1,l_2(X^0_1,Y^1)\rangle+\rho(X^0_1)\langle \alpha^1_2,Y^1\rangle-\langle l_2(X^0_1,X^0_2),\beta^0\rangle+\rho(X^0_2)\langle \alpha^1_1,Y^1\rangle-\langle \alpha^1_1,l_2(X^0_2,Y^1)\rangle\\
 &=&\rho(X^0_1)\langle X^0_2,\beta^0\rangle+\rho(X^0_1)\langle Y^1,\alpha^1_2\rangle.
 \end{eqnarray*}
Thus we have $$\rho(e^0_1)(e^0_2,e^1_1)_+=(e^0_1\circ e^0_2,e^1_1)_++(e^0_2,e^0_1\circ e^1_1)_+.$$
The other cases can be proved similarly, we omit the proof. }

\begin{ex}{\rm
Let $(\g_{-1},\g_0,l_1,l_2,l_3)$ be a Lie $2$-algebra. Denote by $\frkd_0=\g_0\oplus \g^*_{-1}$ and $\frkd_{-1}=\g_{-1}\oplus \g^*_{0}$. Then the  $\LWX$ $2$-algebroid given by Theorem \ref{thm:stL2A} is over a point. By remark \ref{rmk:a point}, we obtain a metric Lie $2$-algebra structure on the graded vector space $\frkd_0\oplus\frkd_{-1}$. The Lie $2$-algebra $(\frkd_{-1},\frkd_0,\partial, [\cdot,\cdot],\Omega)$ is given as follows:
\begin{eqnarray*}
\partial&=&l_1+l^*_1,\\
{[x^0+\alpha^1,y^0+\beta^1]}&=&l_2(x^0,y^0)+{\ad^0}^*_{x^0}\beta^1-{\ad^0}^*_{y^0}\alpha^1,\\
{[x^0+\alpha^1,y^1+\beta^0]}&=&l_2(x^0,y^1)+{\ad^0}^*_{x^0}\beta^0-{\ad^1}^*_{y^1}\alpha^1,\\
\Omega(x^0+\alpha^1,y^0+\beta^1,z^0+\zeta^1)&=&l_3(x^0,y^0,z^0)+{\ad^3}^*_{x^0,y^0}\zeta^1+{\ad^3}^*_{y^0,z^0}\alpha^1+{\ad^3}^*_{z^0,x^0}\beta^1,
\end{eqnarray*}
for all $x^0,y^0,z^0\in\g_{0},$ $x^1,y^1\in\g_{-1},$ $\alpha^1,\beta^1\in\g^*_{-1},$ $\alpha^0,\beta^0\in\g^*_{0}$, where ${\ad^0}^*_{x^0}:\g^*_{-i}\longrightarrow \g^*_{-i}$, ${\ad^1}^*_{x^1}:\g^*_{-1}\longrightarrow \g^*_{0}$ and ${\ad^3}^*_{x^0,y^0}:\g^*_{-1}\longrightarrow\g^*_{0}$ are defined respectively by
\begin{eqnarray*}
\langle{\ad^0}^*_{x^0}\alpha^1,x^1\rangle&=&-\langle\alpha^1,l_2(x^0,x^1)\rangle,\\
\langle{\ad^0}^*_{x^0}\alpha^0,y^0\rangle&=&-\langle\alpha^0,l_2(x^0,y^0)\rangle,\\
\langle{\ad^1}^*_{x^1}\alpha^1,y^0\rangle&=&-\langle\alpha^1,l_2(x^1,y^0)\rangle,\\
\langle{\ad^3}^*_{x^0,y^0}\alpha^1,z^0\rangle&=&-\langle\alpha^1,l_3(x^0,y^0,z^0)\rangle.
\end{eqnarray*}
Thus, this Lie 2-algebra is exactly the semidirect product of the Lie 2-algebra $(\g_{-1},\g_0,l_1,l_2,l_3)$ with its dual $\g_0^*[1]\oplus \g_{-1}^*[1]$ via the coadjoint representation.
}

\emptycomment{
The symmetric bilinear form given by
\begin{eqnarray*}
S(x^0+\alpha^1+x^1+\alpha^0,y^0+\beta^1+y^1+\beta^0)&=&\langle x^0,\beta^0\rangle+\langle\alpha^1,y^1\rangle+\langle x^1,\beta^1\rangle+\langle\alpha^0,y^0\rangle,
\end{eqnarray*}
and the Lie $2$-algebra structure given by
}

\end{ex}

\section{QP-manifolds $T^*[3]A[1]$ and $\LWX$ 2-algebroids}

Let $A$ be a vector bundle over $M$ and $A^*$ its dual bundle. The shifted bundle $A[1]$ is a graded manifold whose fiber space has  degree $-1$. We consider the shifted cotangent bundle $\huaM:=T^*[3]A[1]$. It is a $P$-manifold of degree $3$ over $M$. In this section, we construct a $\LWX$ 2-algebroid   from the degree 3 QP-manifold $T^*[3]A[1]$.

 Denote by $(q^i,\xi^\alpha,\xi_\alpha,p_i)$ a canonical Darboux coordinate on $T^*[3]A[1]$, where $q^i$ is a coordinate on $M$, $\xi^\alpha$ is the fiber coordinate on $A[1]$, $(p_i,\xi_\alpha)$ is the momentum coordinate on $T^*[3]A[1]$ for $(q^i,\xi^\alpha)$. The degrees of variables $(q^i,\xi^\alpha,\xi_\alpha,p_i)$ are respectively $(0,1,2,3)$.
 The degree of the symplectic structure $\omega=dq^idp_i+d\xi^\alpha d\xi_\alpha$ is $3$ and the degree of the corresponding graded Poisson structure is $-3$.
 In the local coordinate, any $Q$-structure $\Theta$ is of the following form:
 \begin{equation}
\Theta={f_1}^i_a(x)p_i\xi^a+{f_2}^{ab}(x)\xi_a\xi_b+\half{f_3}^c_{ab}(x)\xi^a\xi^b\xi_c+\frac{1}{6}{f_4}_{abcd}(x)\xi^a\xi^b\xi^c\xi^d.
\end{equation}
We write $\Theta=\theta_2+\theta_{13}+\theta_4,$ where the substructures are
\begin{eqnarray*}
\theta_2&=&{f_2}^{ab}\xi_a\xi_b,\\
\theta_{13}&=&{f_1}^i_a(x)p_i\xi^a+\half{f_3}^c_{ab}\xi^a\xi^b\xi_c,\\
\theta_4&=&\frac{1}{6}{f_4}_{abcd}\xi^a\xi^b\xi^c\xi^d.
\end{eqnarray*}
The classical master equation $\{\Theta,\Theta\}=0$ is equivalent to the following three identities:
\begin{eqnarray}
\{\theta_{13},\theta_2\}&=&0,\label{eq:master1}\\
\frac{1}{2}\{\theta_{13},\theta_{13}\}+\{\theta_2,\theta_4\}&=&0,\label{eq:master2}\\
\{\theta_{13},\theta_4\}&=&0.\label{eq:master3}
\end{eqnarray}

Define two bundle maps $\partial:A^*\longrightarrow A$ and $\rho:A\longrightarrow TM$  by the following identities respectively:
\begin{eqnarray}
\partial\alpha&=&\{\alpha,\theta_2\},\quad \forall~\alpha\in\Gamma(A^*),\label{eq:Q5}\\
\rho(X)(f)&=&\{f,\{X,\theta_{13}\}\},\quad \forall~X\in\Gamma(A),f\in\CWM\label{eq:Q4}.
\end{eqnarray}
A natural non-degenerate bilinear form $S$ on $A^*\oplus A$ is given by
\begin{equation}
S(X+\alpha,Y+\beta)=\langle X,\beta\rangle+\langle Y,\alpha \rangle,\quad \forall~X,Y\in\Gamma(A),\alpha,\beta\in\Gamma(A^*).\label{eq:naturalsymmetricform}
\end{equation}
Define the   operation $\diamond $ by
\begin{equation}\label{eq:Q1}\left\{\begin{array}{rcll}
X\diamond Y&=&\{Y,\{X,\theta_{13}\}\},&\forall~ X,Y\in\Gamma(A),\\
X\diamond\alpha&=&\{\alpha,\{X,\theta_{13}\}\},&\forall~ X\in\Gamma(A),\alpha\in\Gamma(A^*),\\
\alpha\diamond X&=&-\{X,\{\alpha,\theta_{13}\}\},&\forall~ X\in\Gamma(A),\alpha\in\Gamma(A^*).
\end{array}\right.
\end{equation}
An $A^*$-valued 3-form $\Omega$ is defined by
\begin{eqnarray}
\Omega(X,Y,Z)=\{Z,\{Y,\{X,\theta_4\}\}\},\quad \forall~X,Y,Z\in\Gamma(A).\label{eq:Q7}
\end{eqnarray}
\begin{thm}\label{thm:QPC2A}
  Let $(T^*[3]A[1],\Theta)$ be a $QP$-manifold of degree $3$. Then $( A^*[1],A,\partial, \rho,S,\diamond, \Omega )$ is  a $\LWX$ $2$-algebroid, where $\partial$ is given by \eqref{eq:Q5}, $\rho$ is given by \eqref{eq:Q4}, $S$ is given by   $\eqref{eq:naturalsymmetricform}$, $\diamond$ is given by \eqref{eq:Q1} and  $\Omega$ is given by \eqref{eq:Q7}.
\end{thm}

The proof follows from the following Lemma \ref{lem:skewsymmetry}-\ref{lem:invariant} directly.

\begin{lem}\label{lem:skewsymmetry}
 With the above notations, $e\diamond e=\half \huaD S(e,e)$, where  $S$ is given by   $\eqref{eq:naturalsymmetricform}$ and $\huaD:C^\infty(M)\longrightarrow \Gamma(A^*)$ is given by $\langle\huaD(f),X\rangle=\rho(X)(f)$, in which $\rho$ is given by \eqref{eq:Q4}.
\end{lem}
\pf By \eqref{eq:Q1} and \eqref{eq:delta}, we can deduce that
\begin{eqnarray}
X\diamond Y&=&-Y\diamond X,\\
  X\diamond\alpha+\alpha\diamond X&=&\delta\langle X,\alpha\rangle,\label{eq:Qsym}
\end{eqnarray}
which finishes the proof.\qed

\begin{lem}\label{lem:Leibniz 2}
 With the above notations, $(\Gamma(A^*),\Gamma(A),\partial,\diamond,\Omega)$ is a  Leibniz $2$-algebra, where $\partial$ is given by \eqref{eq:Q5}, $\diamond$ is given by \eqref{eq:Q1} and  $\Omega$ is given by   \eqref{eq:Q7} respectively.
\end{lem}
 \pf
 By  \eqref{eq:master1}, we have $\{\theta_2,\{X,\theta_{13}\}\}=0.$ Thus we have
\begin{eqnarray}
\nonumber\partial(X\diamond\alpha)&=&\{\{\alpha,\{X,\theta_{13}\}\},\theta_2\}=-\{\{\theta_2,\alpha\},\{X,\theta_{13}\}\}-\{\alpha,\{\theta_2,\{X,\theta_{13}\}\}\}\\\label{eq:Qpartial1}
&=&\{\{\alpha,\theta_2\},\{X,\theta_{13}\}\}=X\diamond\partial(\alpha).
\end{eqnarray}
By  \eqref{eq:master1}, we get
\begin{equation}\label{eq:Qrelation}
  \rho\circ\partial=0.
\end{equation}
Then by \eqref{eq:Qsym}, we have
\begin{eqnarray}\label{eq:Qpartial2}
\partial(\alpha\diamond X)&=&\partial(\delta\langle X,\alpha\rangle-X\diamond\alpha)=\partial(\delta\langle X,\alpha\rangle)-X\diamond\partial(\alpha)=\partial(\alpha)\diamond X.
\end{eqnarray}
Similarly, we have
\begin{equation}
  \partial(\alpha)\diamond\beta=\alpha\diamond\partial(\beta)\label{eq:Qpartial3}.
\end{equation}

  By
\eqref{eq:master2} and the following two facts:
\begin{eqnarray*}
\{Z,\{Y,\{X,\{\theta_{13},\theta_{13}\}\}\}\}&=&-2\big({X\diamond(Y\diamond Z)}-(X\diamond Y)\diamond Z- Y\diamond (X\diamond Z)\big),\\
\{Z,\{Y,\{X,\{\theta_{2},\theta_{4}\}\}\}\}&=&\partial\Omega(X,Y,Z),
\end{eqnarray*}
where $X,Y,Z\in\Gamma(A)$,
we have
\begin{eqnarray}
X\diamond (Y\diamond Z)-{(X\diamond Y)\diamond Z}-{Y\diamond(X\diamond Z)}&=&\partial\Omega(X,Y,Z).\label{eq:Qjacobi1}
\end{eqnarray}
\emptycomment{
Eq. \eqref{eq:Qjacobi2} can be proved by
\eqref{eq:master2} and the following two facts:
\begin{eqnarray*}
\{\alpha,\{Y,\{X,\{\theta_{13},\theta_{13}\}\}\}\}&=&-2\big(X\circ(Y\circ\alpha)-(X\circ Y)\circ\alpha-Y\circ(X\circ\alpha)\big),\\
\{\alpha,\{Y,\{X,\{\theta_{2},\theta_{4}\}\}\}\}&=&\Omega(X,Y,\partial\alpha).
\end{eqnarray*}
}
Similarly, we can obtain
\begin{eqnarray}
X\diamond(Y\diamond\alpha)-(X\diamond Y)\diamond\alpha-Y\diamond(X\diamond\alpha)&=&\Omega(X,Y,\partial(\alpha)),\label{eq:Qjacobi2}\\
X\diamond(\alpha\diamond Y)-(X\diamond \alpha)\diamond Y-\alpha\diamond(X\diamond Y)&=&\Omega(X,\partial(\alpha),Y),\label{eq:Qjacobi3}\\
\alpha\diamond(X\diamond Y)-(\alpha\diamond X)\diamond Y-X\diamond(\alpha\diamond Y)&=&\Omega(\partial(\alpha),X,Y).\label{eq:Qjacobi4}
\end{eqnarray}
Finally, expanding $\{W,\{Z,\{Y,\{X,\{\theta_{13},\theta_4\}\}\}\}\}=0$ by the graded Jacobi identity, we have
\begin{eqnarray}\label{eq:Q3form}
&&W\diamond\Omega(X,Y,Z)-X\diamond\Omega(W,Y,Z)+Y\diamond\Omega(W,X,Z)+\Omega(W,X,Y)\diamond Z\nonumber\\
&&-\Omega(W\diamond X,Y,Z)-\Omega(X,W\diamond Y,Z)-\Omega(X,Y,W\diamond Z)\nonumber\\
&&+\Omega(W,X\diamond Y,Z)+\Omega(W,Y,X\diamond Z)-\Omega(W,X, Y\diamond Z)=0.
\end{eqnarray}

By \eqref{eq:Qpartial1}, \eqref{eq:Qpartial2}, \eqref{eq:Qpartial3}, \eqref{eq:Qjacobi1}-\eqref{eq:Q3form}, we deduce that
$(\Gamma(A^*),\Gamma(A),\partial,\diamond,\Omega)$ is a  Leibniz $2$-algebra.\qed

\begin{lem}\label{lem:invariant}
 With the above notations, for all $ \alpha,\beta\in\Gamma(A^*)$, $X,Y,Z,W\in\Gamma(A)$ and $e_1,e_2,e_3\in\Gamma(A)\oplus \Gamma(A^*)$, we have
 \begin{eqnarray}
\label{eq:inv1}\langle\partial\alpha,\beta\rangle&=&\langle\alpha,\partial\beta\rangle, \\
\label{eq:inv2}\rho(e_1)S( e_2,e_3)&=&S( e_1\diamond e_2,e_3)+S(e_2,e_1\diamond e_3), \\
 \label{eq:inv3}S(\Omega(X,Y,Z),W)&=&-S(Z,\Omega(X,Y,W)).
 \end{eqnarray}
\end{lem}
\pf By the Jacobi identity of the graded Poisson bracket $\{\cdot,\cdot\}$, we have
\begin{eqnarray*}
\langle\partial\alpha,\beta\rangle=\{\partial\alpha,\beta\}=\{\{\alpha,\theta_2\},\beta\}=\{\alpha,\{\theta_2,\beta\}\}-\{\theta_2,\{\alpha,\beta\}\}=
-\{\alpha,\partial\beta\}=\{\partial \beta,\alpha\}=\langle\partial \beta,\alpha\rangle.
\end{eqnarray*}

For $X,Y\in\Gamma(A),\alpha\in\Gamma(A^*)$, we have
\begin{eqnarray*}
\{Y,\{\alpha,\{X,\theta_{13}\}\}\}&=&\{\{Y,\alpha\},\{X,\theta_{13}\}\}+\{\alpha,\{Y,\{X,\theta_{13}\}\}\},
\end{eqnarray*}
which implies that
\begin{eqnarray*}
\langle Y,X\diamond\alpha\rangle&=&\rho(X)\langle Y,\alpha\rangle-\langle \alpha,X\diamond Y\rangle.
\end{eqnarray*}
That is $\rho(X)S( Y,\alpha)=S( X\diamond Y,\alpha)+S(Y,X\diamond \alpha).$ Therefore, \eqref{eq:inv2} holds when $e_1,e_2\in \Gamma(A)$ and $e_3\in\Gamma(A^*)$. Similarly, we can show that \eqref{eq:inv2} holds for all the other cases.
\emptycomment{
\begin{eqnarray*}
-\{Y,\{X,\{\alpha,\theta_{13}\}\}\}&=&-\{\{Y,X\},\{\alpha,\theta_{13}\}\}+\{X,\{Y,\{\alpha,\theta_{13}\}\}\}\\
\langle Y,\alpha\circ X\rangle&=&-\langle X,\alpha\circ Y\rangle,
\end{eqnarray*}
This means that $\langle Y,\alpha\circ X\rangle+\langle X,\alpha\circ Y\rangle=0$. }

Finally,  \eqref{eq:inv3} follows from
\begin{eqnarray*}
S(\Omega(X,Y,Z),W)&=&\Poisson{W,\Poisson{Z,\Poisson{Y,\Poisson{X,\theta_4}}}}\\
&=&\Poisson{\Poisson{W,Z},\Poisson{Y,\Poisson{X,\theta_4}}}-\Poisson{Z,\Poisson{W,\Poisson{Y,\Poisson{X,\theta_4}}}}\\
&=&-S(\Omega(X,Y,W),Z).
\end{eqnarray*}
The proof is finished. \qed

\emptycomment{
\begin{pro}
For all $X,Y\in\Gamma(A),\alpha,\beta\in\Gamma(A^*),f\in\CWM$, we have
\begin{eqnarray}
{X\circ fY}&=&f(X\circ Y)+\rho(X)(f)Y,\\
X\circ f\alpha &=&f(X\circ \alpha)+\rho(X)(f)Y,\\
\rho(X\circ Y)&=&[\rho(X),\rho(Y)],
\end{eqnarray}
\end{pro}
 }

\begin{rmk}\label{rmk:QP}
 The P-manifold of degree $3$, $T^*[3]A[1]$, can be viewed as a shifted manifold of $T^*[2]A[1]$, which is a P-manifold of degree $2$. However, in general, a degree $3$ function $\Theta$ on $T^*[2]A[1]$ is not a degree $4$ function on $T^*[3]A[1]$. Thus, there is not a canonical way to obtain a QP-manifold of degree $3$ from a given QP-manifold of degree $2$. Therefore, we can not obtain a $\LWX$ $2$-algebroid from an arbitrary Courant algebroid.
\end{rmk}

\begin{rmk}
 Let us consider the degree $3$ QP-manifold $T^*[3]T[1]M$ where the $Q$-structure is given by $p_i\xi^i$ in local coordinates. On one hand, according to Theorem \ref{thm:QPC2A}, we obtain the $\LWX$ $2$-algebroid $(T^*[1]M,TM,\partial=0,\rho={\id}, S,\diamond,\Omega=0)$ given in Remark \ref{rmk:CC2}. Then according to Theorem \ref{thm:Lie3}, we have a Lie $3$-algebra structure on $C^\infty(M)[2]\oplus \Omega^1(M)[1]\oplus \frkX(M)$. On the other hand, according to \cite{zambon:l-infty}, there is also a Lie $3$-algebra structure on $C^\infty(M)[2]\oplus \Omega^1(M)[1]\oplus (\frkX(M)\oplus \Omega^2(M))$. However, we do not find any connection between the two Lie $3$-algebras.

  Furthermore, if we consider the $Q$-structure  given by $p_i\xi^i+\frac{1}{6}{f_4}_{abcd}\xi^a\xi^b\xi^c\xi^d$, we obtain the $\LWX$ $2$-algebroid $(T^*[1]M,TM,\partial=0,\rho={\id}, S,\diamond,\Omega=H)$ given in Example \ref{ex:4-form}.
\end{rmk}

\section{The $\LWX$ $2$-algebroid associated to a split Lie 2-bialgebroid}
In this section, we introduce the notion of a split Lie 2-bialgebroid and show that there is a $\LWX$ 2-algebroid structure on $\huaA\oplus\huaA^*[1]$ associated to any split  Lie 2-bialgebroid $(\huaA,\huaA^*[1])$.

Now assume that there is a  split Lie $2$-algebroid structure on the dual bundle $\huaA^*[1]=A^*_0[1]\oplus A^*_{-1}[1]$. Since $T^*[3]((A_0\oplus A^*_{-1})^*[1])[1]$, $T^*[3](A_0\oplus A^*_{-1})[1]$ and $T^*[3](A_0\oplus A^*_{-1})^*[2]$ are naturally isomorphic, by Theorem \ref{thm:L2A-MST}, the dual split Lie $2$-algebroid $(\huaA^*[1],\frkl_1,\frkl_2,\frkl_3,\frka)$ gives rise to a degree $4$ function $\gamma$ on $T^*[3](A_0\oplus A^*_{-1})[1]$ satisfying $\Poisson{\gamma,\gamma}=0.$ It is given in local coordinates $(x^i,\xi^j,\theta_k, p_i,\xi_j,\theta^k)$ by
\begin{equation}\label{eq:gamma}
\gamma={\gamma_1}^{ij}(x)p_j\theta_i+{\gamma_2}^{j}_i(x)\xi_j\theta^i+\half{\gamma_3}_k^{ij}(x)\theta^k\theta_i\theta_j+ {\gamma_4}_{k}^{ij}\xi_i\theta_j\xi^k+\frac{1}{6}{\gamma_5}^{ijk}_l(x)\xi^l\theta_i\theta_j\theta_k.
\end{equation}
We will also write $\gamma=\gamma_2+\gamma_{134}+\gamma_5.$

\begin{defi}
Let $(\huaA,l_1,l_2,l_3,a)$ be a split Lie $2$-algebroid with the structure function $\mu$ given by \eqref{eq:L2Amu} and $(\huaA^*[1],\frkl_1,\frkl_2,\frkl_3,\frka)$ a split Lie $2$-algebroid with the structure function $\gamma$ given by \eqref{eq:gamma}. The pair $(\huaA,\huaA^*[1])$ is called a {\bf split Lie $2$-bialgebroid} if $\gamma_2=\mu_2$ and
\begin{equation}
\Poisson{\mu+\gamma-\mu_2,\mu+\gamma-\mu_2}=0,
\end{equation}
where $\Poisson{\cdot,\cdot}$ is the graded Poisson bracket corresponding to the symplectic structure $\omega=dx^idp_i+d\xi^j d\xi_j+d\theta_kd\theta^k$ on  $T^*[3](A_0\oplus A^*_{-1})[1]$.
\end{defi}
Denote a split Lie $2$-bialgebroid by $(\huaA,\huaA^*[1])$.

We denote by $\huaL^0,\huaL^1,\huaL^3,\delta_*$ the operations for the dual split Lie $2$-algebroid $(\huaA^*[1],\frkl_1,\frkl_2,\frkl_3,\frka)$ corresponding to the operations $L^0,L^1,L^3,\delta$ for the split Lie $2$-algebroid $(\huaA,l_1,l_2,l_3,a)$.

Now we assume that $(\huaA,l_1,l_2,l_3,a)$  and $(\huaA^*[1],\frkl_1,\frkl_2,\frkl_3,\frka)$ are split Lie $2$-algebroids.
Let $E_0=A_0\oplus A^*_{-1}$, $E_{-1}=A_{-1}\oplus A^*_{0}$ and $\huaE=E_0\oplus E_{-1}$.

Let $\partial:E_{-1}\longrightarrow E_0$ and $\rho:E_0\longrightarrow TM$ be bundle maps defined by
\begin{eqnarray}
\label{eq:parbi}\partial(X^1+\alpha^0)&=&l_1(X^1)+\frkl_1(\alpha^0),\label{eq:L2Bpartial}\\
\rho(X^0+\alpha^1)&=&a(X^0)+\frka(\alpha^1)\label{eq:L2Banchor}.
\end{eqnarray}

On $\Gamma(\huaE)$, we introduce the operation $\diamond$ by
\begin{equation}\label{eq:L2Bbracket0}\left\{\begin{array}{rcl}
(X^0+\alpha^1)\diamond(Y^0+\beta^1)&=&l_2(X^0,Y^0)+L^0_{X^0}\beta^1-L^0_{Y^0}\alpha^1 +\frkl_2(\alpha^1,\beta^1)+\huaL^0_{\alpha^1}Y^0-\huaL^0_{\beta^1}X^0,\\
(X^0+\alpha^1)\diamond(X^1+\alpha^0)&=&l_2(X^0,X^1)+L^0_{X^0}\alpha^0+\iota_{X^1}\delta(\alpha^1) +\frkl_2(\alpha^1,\alpha^0)+\huaL^0_{\alpha^1}X^1+\iota_{\alpha^0}\delta_*(X^0), \\
(X^1+\alpha^0)\diamond(X^0+\alpha^1)&=&l_2(X^1,X^0)+L^1_{X^1}\alpha^1-\iota_{X^0}\delta(\alpha^0) +\frkl_2(\alpha^0,\alpha^1)+\huaL^1_{\alpha^0}X^0-\iota_{\alpha^1}\delta_*(X^1).
\end{array}\right.
\end{equation}

An $E_{-1}$-valued $3$-form $\Omega$ is defined by
\begin{eqnarray}\label{eq:L2B3-form}
\nonumber\Omega(X^0+\alpha^1,Y^0+\beta^1,Z^0+\zeta^1)
&=&l_3(X^0,Y^0,Z^0)+L^3_{X^0,Y^0}\zeta^1+L^3_{Y^0,Z^0}\alpha^1+L^3_{Z^0,X^0}\beta^1\\
&&+\frkl_3(\alpha^1,\beta^1,\zeta^1)+\huaL^3_{\alpha^1,\beta^1}Z^0+\huaL^3_{\beta^1,\zeta^1}X^0+\huaL^3_{\zeta^1,\alpha^1}Y^0,
\end{eqnarray}
where $X^0,Y^0,Z^0\in\Gamma(A_0),\alpha^1,\beta^1,\zeta^1\in\Gamma(A_{-1}^*).$

\begin{thm}\label{thm:Lie2biC2}
Let $(\huaA,\huaA^*[1])$ be a split Lie $2$-bialgebroid. Then $(E_{-1},E_0,\partial,\rho,(\cdot,\cdot)_+,\diamond,\Omega)$ is a $\LWX$ $2$-algebroid, where  $E_0=A_0\oplus A^*_{-1}$, $E_{-1}=A_{-1}\oplus A^*_{0}$, $\partial$ is given by \eqref{eq:L2Bpartial}, $\rho$ is given by \eqref{eq:L2Banchor}, $(\cdot,\cdot)_+$ is given by \eqref{eq:naturalsymform}, $\diamond$ is given by \eqref{eq:L2Bbracket0} and $\Omega$ is given by \eqref{eq:L2B3-form}.
\end{thm}
\pf Since $\mu+\gamma-\mu_2$ is a degree 4 function on $T^*[3]E_0[1]$ satisfying $\Poisson{\mu+\gamma-\mu_2,\mu+\gamma-\mu_2}=0$, by Theorem \ref{thm:QPC2A}, there is a $\LWX$ 2-algebroid  defined by $\mu+\gamma-\mu_2$ through derived brackets.  It is straightforward to deduce that   \eqref{eq:parbi}-\eqref{eq:L2B3-form} are exactly the one obtained through derived brackets. The proof is finished.\qed

\end{document}